\documentclass[11pt,a4paper]{amsart}
\usepackage{amsmath,amsthm,amssymb}
\usepackage[all]{xy}
\DeclareMathOperator{\id}{id}

\DeclareMathOperator{\im}{im}

\DeclareMathOperator{\pol}{pol}
\DeclareMathOperator{\red}{red}
\DeclareMathOperator{\Rot}{rot}

\DeclareMathOperator{\Glob}{G}
\DeclareMathOperator{\Loc}{{H}}

\DeclareMathOperator{\Amb}{Amb}
\DeclareMathOperator{\PreAmb}{PreAmb}
\DeclareMathOperator{\CoAmb}{CoAmb}
\DeclareMathOperator{\PreCoAmb}{PreCoAmb}
\DeclareMathOperator{\vFQ}{vFQ}
\DeclareMathOperator{\sFQ}{sFQ}
\DeclareMathOperator{\psFQ}{psFQ}

\DeclareMathOperator{\CH}{CH}
\DeclareMathOperator{\CS}{CS}
\DeclareMathOperator{\CE}{CE}

\newcommand{\covarkappa}{{\bar\varkappa}}
\newcommand{\colambda}{{\bar\lambda}}
\newcommand{\copi}{{\bar\pi}}
\newcommand{\coeta}{{\bar\eta}}
\newcommand{\opp}{\mathrm{opp}}
\swapnumbers
\theoremstyle{definition}
\newtheorem{point}{}[section]
\newtheorem{remark}[point]{Remark}

\newtheorem{example}[point]{Example}
\newtheorem{conven}[point]{Convention}
\theoremstyle{plain}

\newtheorem{lemma}[point]{Lemma}
\newtheorem{theorem}[point]{Theorem}

\newtheorem{cor}[point]{Corollary}
\newcommand{\marginextend}[1]{ \addtolength{\oddsidemargin}{-#1}  \addtolength{\evensidemargin}{-#1}\addtolength{\textwidth}{#1}\addtolength{\textwidth}{#1}}
\newcommand{\updownextend}[1]{ \addtolength{\topmargin}{-#1}  \addtolength{\textheight}{#1}
\addtolength{\textheight}{#1}}
\DeclareMathOperator{\Id}{Id}

\DeclareMathOperator{\Mult}{Mult}
\marginextend{1.25cm}
\updownextend{0.25cm}\allowdisplaybreaks[4]
\title[FQ operations I.
Basic invariance properties]{Fermionic quantum operations: a computational framework I.
Basic invariance properties}
\author{Gyula Lakos}
\date{\today}
\email{lakos@cs.elte.hu}
\thanks{The author would like to thank  Bal\'azs Csik\'os.}
\address{Department of Geometry, E\"otv\"os Lor\'and University, P\'azm\'any P\'eter s.~1/C,  Budapest, H--1117, Hungary}
\keywords{Clifford systems, formal power series, calculus of non-commuting operators}
\subjclass[2010]{Primary: 47L99. Secondary: 15B99, 81Q99.}
\begin{document}
\maketitle
\vspace{-1cm}
\begin{abstract}
The objective of this series of papers is to recover information regarding the behaviour of FQ operations
in the case $n=2$, and  FQ conform-operations in the case $n=3$.
In this first part we study how the basic invariance properties
of FQ operations ($n=2$) are
reflected in their formal power series expansions.
\end{abstract}
\section*{Introduction}
Dealing with functions of several noncommuting operators has many approaches, like
holomorphic calculus around the joint spectrum \cite{Tay}, operator ordering \cite{NSS},
using Clifford variables \cite{BDS}, integrated functional calculi
\cite{Jef}, advanced perturbation techniques and rational noncommutative functions \cite{KVV} (also see references therein).
Here we intend to inves\-ti\-gate a very specific one, basically considering functions defined on perturbations of Clifford systems:
FQ (conform-)operations were introduced by the author
in \cite{L1} as a kind of noncommutative linear algebra and/or functional calculus.
These operation have both analytic and algebraic aspects, here we center on some effective computational techniques regarding them,
primarily from algebraic viewpoint.
The objective of this series of papers is to recover information regarding the behaviour of FQ operations
in the case $n=2$, and  FQ conform-operations in the case $n=3$.
In this first part we study how the basic invariance properties
of FQ operations ($n=2$) are reflected in their formal power series expansions.
We concentrate on formal FQ operations, and, for the sake of simplicity, only on ones with good sign-linear properties.
On the other hand, although our ultimate objective is the study of natural FQ operations, we
will typically consider expansions around a fixed Clifford system $(Q_1,Q_2)$
without the a priori assumption of naturality / conjugation invariance.

\section{FQ operations: expansions and bases}
\begin{point}
Following \cite{L1}, formal FQ operations can be represented as formal power series expansions around Clifford systems.
In the case $n=2$, we have a Clifford system $(Q_1,Q_2)$, and then the FQ operation is considered to make sense for
the pair $(A_1,A_2)=(Q_1+R_1,Q_2+R_2)$, where $R_1$ and $R_2$ are imagined as infinitesimal or formal variables.

If $Q$ is a skew-involution,  then we use notation
\[R{}^{0}_{Q}:=\frac12 \left(R + Q^{-1}RQ\right),
\qquad
R{}^{1}_{Q}:=\frac12 \left(R - Q^{-1}RQ\right);\]
and
\[(R/Q)^{(\iota_1,\iota_2)}_j:=(R_jQ^{-1}_j){}_{Q_1}^{\iota_1}{}_{Q_2}^{\iota_2}.\]

Now, a scalar FQ operation with property \texttt{(sL)} can be written as
\begin{equation} \Upsilon(A_1,A_2)=f_0( (R/Q)  ),\label{f:scalar}\end{equation}
where $f$ is a formal real noncommutative power series in the $8$ terms $(R/Q)^{(\iota_1,\iota_2)}_j$.

Similarly, a vectorial FQ operation with property \texttt{(vL)} can be written as
\begin{equation} \Psi(A_1,A_2)=(f_1( (R/Q)  )Q_1,f_2( (R/Q)  )Q_2);\label{f:vec}\end{equation}
and a  pseudoscalar FQ operation with property \texttt{(psL)} can be written as
\begin{equation} \Phi(A_1,A_2)=f_{12}( (R/Q)  )Q_1Q_2.\label{f:pseudo}\end{equation}

In order to arrive to natural FQ operations, one needs one more assumption, naturality, which is equivalent to uniform analiticy
(as opposed to the pointed expansions above), which is also equivalent to conjugation invariance.
The reader is advised to look up the discussion in \cite{L1}; although  we will also address this question later.

Nevertheless, in this paper, our starting point will be simply considering expansions as in  \eqref{f:scalar}, \eqref{f:vec}, \eqref{f:pseudo}.
Base-point invariance will be required ultimately, but our basic objects will be the pointed expansions as above.
So, strictly speaking, we should always use the notation
$\Xi_{(Q_1,Q_2)}(A_1,A_2)$
in order to indicate that the expansion is taken around the Clifford system $(Q_1,Q_2)$,
but, in general, we will simply write $\Xi(A_1,A_2)$.
(Later, we will often consider the modification
$\Xi^{\mathrm{ext}}(A_1,A_2)=\Xi_{\underline{\mathcal O}^{\mathrm{Sy}}(A_1,A_2)}(A_1,A_2)$,
which will be our method of choice for obtaining natural operations in this paper.)
\end{point}
\begin{conven}
According to this, in what follows, an FQ operation will mean a
pointed expansion as above $(n=2)$, either scalar, vectorial, or pseudoscalar; the sign-linear property
\texttt{(sL)}/\texttt{(vL)}/\texttt{(psL)} built into the expansion.

\end{conven}
\begin{point}
(a) In order to deal with the expansions more effectively, we will use the notation
\[r_1=(R/Q)^{(00)}_1,\quad r_2=(R/Q)^{(01)}_1,\quad r_3=(R/Q)^{(10)}_1,\quad r_4=(R/Q)^{(11)}_1,\]
\[r_5=(R/Q)^{(00)}_2,\quad r_6=(R/Q)^{(01)}_2,\quad r_7=(R/Q)^{(10)}_2,\quad r_8=(R/Q)^{(11)}_2.\]
This means that the free algebra $\mathfrak F_2$ generated by the Clifford elements $Q_1,Q_2$ and formal variables $R_1,R_2$
can also be interpreted as generated by the Clifford elements $Q_1,Q_2$ and formal variables $r_1,\ldots,r_8$.
In the latter case there are more variables but with (anti)commutation rules with respect to $Q_1,Q_2$.
We call the infinitesimal base $r_1,\ldots,r_8$ as the split base $(R/Q)_{\mathrm{split}}$.

(b) Other choice is the mixed basis $(R/Q)_{\mathrm{mix}}$ given by  $\hat{r}_1,\ldots,\hat{r}_8$, where
\begin{equation}
\begin{bmatrix}\hat{r}_1\\\hat{r}_2\\\hat{r}_3\\\hat{r}_4\\\hat{r}_5\\\hat{r}_6\\\hat{r}_7\\\hat{r}_8\end{bmatrix}
= \frac12\begin{bmatrix}
0&1&0&0&0&0&-1&0\\
0&1&0&0&0&0&1&0\\
1&0&0&0&1&0&0&0\\
1&0&0&0&-1&0&0&0\\
0&0&0&1&0&0&0&-1\\
0&0&0&1&0&0&0&1\\
0&0&1&0&0&1&0&0\\
0&0&1&0&0&-1&0&0\end{bmatrix}
\begin{bmatrix}r_1\\r_2\\r_3\\r_4\\r_5\\r_6\\r_7\\r_8\end{bmatrix}.
\label{eq:atter1}\end{equation}
Here the basis elements are mixed from $R_1,R_2$ along ``characters'':
\begin{equation}
\begin{bmatrix}\hat{r}_1\\\hat{r}_2\\\hat{r}_3\\\hat{r}_4\\\hat{r}_5\\\hat{r}_6\\\hat{r}_7\\\hat{r}_8\end{bmatrix}
= \frac18\begin{bmatrix}
1&1&-1&-1&-1&-1&1&1\\
1&1&-1&-1&1&1&-1&-1\\
1&1&1&1&1&1&1&1\\
1&1&1&1&-1&-1&-1&-1\\
1&-1&-1&1&-1&1&1&-1\\
1&-1&-1&1&1&-1&-1&1\\
1&-1&1&-1&1&-1&1&-1\\
1&-1&1&-1&-1&1&-1&1\end{bmatrix}
\begin{bmatrix}R_1Q^{-1}_1\\Q^{-1}_1 R_1\\Q_2R_1Q_1Q_2\\Q_2Q_1R_1Q_2\\
R_2Q^{-1}_2\\ Q^{-1}_2 R_2\\Q_1R_2Q_2Q_1\\Q_1Q_2R_2Q_1 \end{bmatrix}.
\label{eq:charac}
\end{equation}
(The order of elements in the basis is somewhat arbitrary).
The commutation rules for this basis are given by
\[Q_1\begin{bmatrix}\hat{r}_1&\hat{r}_2&\hat{r}_3&\hat{r}_4&\hat{r}_5&\hat{r}_6&\hat{r}_7&\hat{r}_8\end{bmatrix}Q_1^{-1}=
\begin{bmatrix}\hat{r}_2&\hat{r}_1&\hat{r}_3&\hat{r}_4&-\hat{r}_5&-\hat{r}_6&-\hat{r}_8&-\hat{r}_7\end{bmatrix},\]
\[Q_2\begin{bmatrix}\hat{r}_1&\hat{r}_2&\hat{r}_3&\hat{r}_4&\hat{r}_5&\hat{r}_6&\hat{r}_7&\hat{r}_8\end{bmatrix}Q_2^{-1}=
\begin{bmatrix}-\hat{r}_2&-\hat{r}_1&\hat{r}_3&\hat{r}_4&-\hat{r}_5&-\hat{r}_6&\hat{r}_8&\hat{r}_7\end{bmatrix},\]
\[Q_1Q_2\begin{bmatrix}\hat{r}_1&\!\!\!\hat{r}_2&\!\!\hat{r}_3&\!\!\hat{r}_4&\!\!\hat{r}_5&\!\!\hat{r}_6&\!\!\hat{r}_7&\!\!\hat{r}_8\end{bmatrix}(Q_1Q_2)^{-1}=
\begin{bmatrix}-\hat{r}_1&\!-\hat{r}_2&\!\hat{r}_3&\!\hat{r}_4&\!\hat{r}_5&\!\hat{r}_6&\!-\hat{r}_7&\!-\hat{r}_8\end{bmatrix}.\]

(c) A slight variation is  the circular basis $(R/Q)_{\mathrm{circ}}$ given by $\tilde{r}_1,\ldots,\tilde{r}_8$, where
\begin{equation}
\begin{bmatrix}\tilde{r}_1\\\tilde{r}_2\\\tilde{r}_3\\\tilde{r}_4\\\tilde{r}_5\\\tilde{r}_6\\\tilde{r}_7\\\tilde{r}_8\end{bmatrix}
=\begin{bmatrix}
1&0&0&0&0&0&0&0\\
0&1&0&0&0&0&0&0\\
0&0&1&0&0&0&0&0\\
0&0&0&1&1&0&0&0\\
0&0&0&1&-1&0&0&0\\
0&0&0&0&0&1&0&0\\
0&0&0&0&0&0&1&0\\
0&0&0&0&0&0&0&1
\end{bmatrix}
\begin{bmatrix}\hat{r}_1\\\hat{r}_2\\\hat{r}_3\\\hat{r}_4\\\hat{r}_5\\\hat{r}_6\\\hat{r}_7\\\hat{r}_8\end{bmatrix}.
\label{eq:atter2}\end{equation}
The reason for the name will be clear later. The commutation rules for this basis are
\[Q_1\begin{bmatrix}\tilde{r}_1&\tilde{r}_2&\tilde{r}_3&\tilde{r}_4&\tilde{r}_5&\tilde{r}_6&\tilde{r}_7&\tilde{r}_8\end{bmatrix}Q_1^{-1}=
\begin{bmatrix}\tilde{r}_2&\tilde{r}_1&\tilde{r}_3&\tilde{r}_5&\tilde{r}_4&-\tilde{r}_6&-\tilde{r}_8&-\tilde{r}_7\end{bmatrix},\]
\[Q_2\begin{bmatrix}\tilde{r}_1&\tilde{r}_2&\tilde{r}_3&\tilde{r}_4&\tilde{r}_5&\tilde{r}_6&\tilde{r}_7&\tilde{r}_8\end{bmatrix}Q_2^{-1}=
\begin{bmatrix}-\tilde{r}_2&-\tilde{r}_1&\tilde{r}_3&\tilde{r}_5&\tilde{r}_4&-\tilde{r}_6&\tilde{r}_8&\tilde{r}_7\end{bmatrix},\]
\[Q_1Q_2\begin{bmatrix}\tilde{r}_1&\!\!\tilde{r}_2&\!\!\tilde{r}_3&\!\!\tilde{r}_4&\!\!\tilde{r}_5&\!\!\tilde{r}_6&\!\!\tilde{r}_7&\!\!\tilde{r}_8
\end{bmatrix}(Q_1Q_2)^{-1}=
\begin{bmatrix}-\tilde{r}_1&\!-\tilde{r}_2&\!\tilde{r}_3&\!\tilde{r}_4&{}\!\tilde{r}_5&\!\tilde{r}_6&\!-\tilde{r}_7&\!-\tilde{r}_8\end{bmatrix}.\]
\end{point}
\begin{point}
If  $f$ is one of $f_0,f_1,f_2,f_{12}$, then we use the notation
\[f_s((R/Q))\equiv f_s(r_1,\ldots,r_8 )=\sum_{r=0}^\infty\sum_{k_1,\ldots,k_r\in\{1,\ldots,8\}} p^{[s]}_{k_1,\ldots,k_r} r_{k_1}\ldots r_{k_r}.\]

I.~e., when we take noncommutative power series in $(R/Q)$, then (somewhat loosely) it will be understood as a noncommutative power series
in the (order of the) variables of $(R/Q)_{\mathrm{split}}$.
These expressions can be realized by other power series $\hat f_s$, $\tilde f_s$ in the mixed and circular bases.
So, then $f((R/Q)_{\mathrm{split}})=\hat f((R/Q)_{\mathrm{mix}})=\tilde f((R/Q)_{\mathrm{circ}})$.
The coefficients of
$f,\hat f,\tilde f$ can be expressed from each other, using (the inverses of the) matrices from \eqref{eq:atter1}  and \eqref{eq:atter2}.
In particular,
\[\tilde  p^{[s]}_4=\frac12(\hat  p^{[s]}_4+\hat  p^{[s]}_5)\qquad\text{and}\qquad \tilde  p^{[s]}_5=\frac12(\hat  p^{[s]}_4-\hat  p^{[s]}_5).\]

For practical purposes, we will also use the notation $(-1)^{[s]}$, where
$(-1)^{[0]}=(-1)^{[1]}=(-1)^{[12]}=1$ and $(-1)^{[2]}=-1$.
\end{point}
\begin{example}
(a) The expansion of the constant $1$ operation (as a scalar operation), in the split basis, is given by
\[{\mathbf P}^{[1]}_0=[1],\]
and all other coefficients, in expansion orders $r\geq1$, are $0$. The same applies in the mixed and circular bases, too.

(b) The expansion of the identity operation $\Id$ (as a vectorial operation), in the mixed basis, is given by
\[\hat{\mathbf P}^{[1]}_0=[1],\qquad\hat{\mathbf P}^{[1]}_1=\left[ \begin {array}{cccccccc} 1&1&1&1&1&1&1&1\end {array} \right],\]
\[\hat{\mathbf P}^{[2]}_0=[1],\qquad\hat{\mathbf P}^{[2]}_1= \left[ \begin {array}{cccccccc} -1&1&1&-1&-1&1&1&-1\end {array} \right] ,\]
and all other coefficients, in expansion orders $r\geq2$, are $0$. In the circular basis it is given by
\[\tilde{\mathbf P}^{[1]}_0=[1],\qquad\tilde{\mathbf P}^{[1]}_1=\left[ \begin {array}{cccccccc} 1&1&1&1&0&1&1&1\end {array} \right],\]
\[\tilde{\mathbf P}^{[2]}_0=[1],\quad\tilde{\mathbf P}^{[2]}_1= \left[ \begin {array}{cccccccc} -1&1&1&-1&0&1&1&-1\end {array} \right] ,\]
and all other coefficients, in expansion orders $r\geq2$, are $0$.

(c) The expansion of the pseudodeterminant operation $\mathcal D(A_1,A_2)=\frac12[A_1,A_2]$ (as a pseudoscalar operation), in the mixed basis, is given by
\[\hat{\mathbf P}^{[12]}_0=[1],\qquad\hat{\mathbf P}^{[12]}_1=\left[ \begin {array}{cccccccc} 0&0&2&0&0&0&2&0\end {array} \right] ,\]
\[\hat{\mathbf P}^{[12]}_2=\left[ \begin {array}{cccccccc} 1&0&0&-1&1&0&0&-1
\\ 0&-1&1&0&0&-1&1&0\\ 0&-1&1&0&0&
-1&1&0\\ 1&0&0&-1&1&0&0&-1\\ 1&0&0
&-1&1&0&0&-1\\ 0&-1&1&0&0&-1&1&0
\\ 0&-1&1&0&0&-1&1&0\\ 1&0&0&-1&1&0
&0&-1\end {array} \right];\]
the other coefficients, in expansion orders $r\geq3$, are $0$.
\end{example}
\begin{point}
Consider the (conform-)orthogonalization procedures
$\mathcal O^{\mathrm{GS}}$, $\mathcal O^{\mathrm{Sy}}$, $\mathcal O^{\mathrm{fGS}}$, $\mathcal O^{\mathrm{fSy}}$
from \cite{L1}. Then one can easily see the following:

$(Q_1,Q_2)$ is the Gram-Schmidt orthogonalization of $(A_1,A_2)=(Q_1+R_1,Q_2+R_2)$ if
and only if $r_3=r_4=r_6=0$ in the split base.

$(Q_1,Q_2)$ is the symmetric orthogonalization of $(A_1,A_2)=(Q_1+R_1,Q_2+R_2)$ if
and only if $\hat{r}_6=\hat{r}_7=\hat{r}_8=0$ in the mixed base; and the same applies to the circular base.

$(Q_1,Q_2)$ is the Gram-Schmidt conform-orthogonalization of $(A_1,A_2)=(Q_1+R_1,Q_2+R_2)$ if
and only if  $r_1=r_2=r_3=r_4=r_6=r_8=0$ in the split base.

$(Q_1,Q_2)$ is the symmetric conform-orthogonalization of $(A_1,A_2)=(Q_1+R_1,Q_2+R_2)$ if
and only if $\hat{r}_1=\hat{r}_2=\hat{r}_3=\hat{r}_6=\hat{r}_7=\hat{r}_8=0$ in the mixed base; and the same applies to the circular base.
\end{point}
\begin{point} The simplest invariance property for FQ operations is

$\texttt{(sC)}$ / $\texttt{(vC)}$ / $\texttt{(psC)}$  Clifford conservativity:
It means $\Xi(Q_1,Q_2)=1$ /  $\Xi(Q_1,Q_2)=(Q_1,Q_2)$ / $\Xi(Q_1,Q_2)=Q_1Q_2$ respectively.

In terms of the expansions, it can be expressed by
$p^{[0]}=1$ / $p^{[1]}=p^{[2]}=1$  / $p^{[12]}=1$  respectively  (and the same in other bases),
i.~e., the leading coefficients of the power series expansion are $1$.
Together with naturality / conjugation invariance, this implies that $\Xi$ yields the expected  simple results on any
Clifford system.
We will often (but not always) assume this property.
\end{point}

\section{Principal invariance properties}
Here we study the most basic invariance properties.
The discussion might seem to be a bit redundant
regarding the use of various bases, but it is useful to see
how these bases are different from each other.

\begin{point} \textbf{Conjugation invariance \texttt{(Nat)}.}\label{po:conj}
This is equivalent to naturality, and base Clifford system invariance.
I.~e.~in the expansions \eqref{f:scalar}, \eqref{f:vec}, \eqref{f:pseudo}
one can use Gram-Schmidt orthogonalization or the symmetric orthogonalization as base system.

In terms of the split basis and the Gram-Schmidt orthogonalization, this means that in the expansions \eqref{f:scalar}, \eqref{f:vec}, \eqref{f:pseudo},
the formal power series $f(r_1,\ldots,r_8)$ can be
reconstructed uniquely from $f(r_1,r_2,0,0,r_5,0,r_7,r_8)$, and that any formal power series
$g(r_1,r_2,r_5,r_7,r_8)$ can be prescribed for the latter.
Indeed, if $A=Q+R$, $Q^{\mathrm{GS}}=\underline{\mathcal O}^{\mathrm{GS}}(A)$ and $A=Q^{\mathrm{GS}}+R^{\mathrm{GS}}$,
\[(R/Q)_{\mathrm{split}}=\begin{bmatrix}{r}_1&{r}_2&{r}_3&{r}_4&{r}_5&{r}_6&{r}_7&{r}_8\end{bmatrix},\]
then one case see that
\begin{multline}
(R^{\mathrm{GS}}/Q^{\mathrm{GS}})_{\mathrm{split}}=
[\begin{matrix}{r}_1+O(R/Q)^2&{r}_2+O(R/Q)^2&0&\ldots\end{matrix}\label{e:split}\\
\begin{matrix}\ldots&0&{r}_5+O(R/Q)^2&0&{r}_7+O(R/Q)^2&{r}_8-{r}_4+O(R/Q)^2\end{matrix}].\end{multline}
Hence, when  \eqref{f:scalar}, \eqref{f:vec}, \eqref{f:pseudo} are expanded around $Q^{\mathrm{GS}}$ (which we can compute)
only the ${1,2,5,7,8}$-coefficients count.
This argument is worked out in \cite{L1} in the vectorial case, up to order $2$.

In terms of the mixed basis and symmetric orthogonalization, this means that in the expansions
the formal power series $f(\hat{r}_1,\ldots,\hat{r}_8)$ can be
reconstructed from the restrictions $f(\hat{r}_1,\hat{r}_2,\hat{r}_3,\hat{r}_4,\hat{r}_5,0,0,0)$, and also that any formal power series
$g(\hat{r}_1,\hat{r}_2,\hat{r}_3,\hat{r}_4,\hat{r}_5)$ can be prescribed for the latter.
Using the mixed basis, if $A=Q+R$, $Q^{\mathrm{Sy}}=\underline{\mathcal O}^{\mathrm{Sy}}(A)$ and $A=Q^{\mathrm{Sy}}+R^{\mathrm{Sy}}$,
\[(R/Q)_{\mathrm{mix}}=\begin{bmatrix}\hat{r}_1&\hat{r}_2&\hat{r}_3&\hat{r}_4&\hat{r}_5&\hat{r}_6&\hat{r}_7&\hat{r}_8\end{bmatrix},\]
then one case see that
\begin{multline}(R^{\mathrm{Sy}}/Q^{\mathrm{Sy}})_{\mathrm{mix}}=
[\begin{matrix}\hat{r}_1+O(R/Q)^{\geq2}&\hat{r}_2+O(R/Q)^{\geq2}&\hat{r}_3+O(R/Q)^{\geq2}\end{matrix}\ldots\\
\ldots\begin{matrix}&\hat{r}_4+O(R/Q)^{\geq2}&\hat{r}_5+O(R/Q)^{\geq2}&0&0&0\end{matrix}].\label{e:mix}\end{multline}
Hence, when  \eqref{f:scalar}, \eqref{f:vec}, \eqref{f:pseudo} are expanded around $Q^{\mathrm{Sy}}$ (which we can compute)
only the ${1,2,3,4,5}$-coefficients count.
The very same argument applies to the circular basis.

This means that for (pseudo)scalar operations we have $5^r$ coefficients to be chosen freely in the $r$th perturbation level
(the pure $\{1,2,5,7,8\}$-terms in the split basis, and the pure $\{1,2,3,4,5\}$-terms in the mixed/circular bases),
and $2\times 5^r$ coefficients in the vectorial case; and all the other coefficients can be expressed from them explicitly.
The actual reductions/extensions are to be found from
\[\Xi_{(Q_1,Q_2)}(Q_1+R_1,Q_2+R_2)=\Xi_{\underline{\mathcal O}^\omega(Q_1+R_1,Q_2+R_2)}(Q_1+R_1,Q_2+R_2),\]
but they are not particularly simple. (Although, later, we will show explicit recursion relations  for them by direct methods.)

There is a small advantage using the mixed/circular bases to the split basis.
This is as follows:
When we consider an expansion according to \eqref{e:split},
and we eliminate the $\{3,4,6\}$-indices,
we see that in the result the coefficient of $r_4$ gets a contribution from the coefficient of $r_8$.
This is reflected in the natural coefficient rule
$p_4^{[1]}=p^{[1]}-p_8^{[1]}$
(in the split basis), cf. \cite{L1}. On the other hand, in the mixed/circular bases the coefficients of the eliminated indices $\{6,7,8\}$
do not receive contribution from the same order; hence, during  the reduction process the coefficient indices  $\{6,7,8\}$ ``decay'' to
a lower order:

\end{point}
\begin{example} Similarly, to the example in \cite{L1},
let $\Psi$ be a formal FQ operation with \texttt{(vL)}, $n=2$; but now we use the circular base. Then
\[\Psi(Q+R)_s=\left(\tilde p^{[s]}+\sum_{i=1}^8 \tilde p^{[s]}_i\tilde r_i+\sum_{i,j=1}^8 \tilde p^{[s]}_{ij}\tilde r_i\tilde r_j+O((R/Q)^{\geq3}) \right)Q_s,\qquad s\in\{1,2\}. \tag{F${}_2$}\label{U}\]

We collect (some of) the coefficients into the scalar matrices $\mathbf P^{[s]}_0=[ \tilde p^{[s]}]$, the row matrices
$\tilde{\mathbf P}^{[s]}_1=[ \tilde p^{[s]}_i]_{i=1}^8$, and the square matrices $\tilde{\mathbf P}^{[s]}_2=[ \tilde p^{[s]}_{ij}]_{i=1}^8{}_{j=1}^8$.

If $\Psi=\underline{\mathcal O}^{\mathrm{Sy}}$, then one can check that
$\tilde{\mathbf P}^{[1]}_0=[1]$, $ \tilde{\mathbf P}^{[2]}_0=[1]$, and
\[\tilde{\mathbf P}^{[1]}_1=\left[ \begin {array}{cccccccc} 0&0&0&0&0&1&1&1\end {array} \right],\]
\[\tilde{\mathbf P}^{[2]}_1= \left[ \begin {array}{cccccccc} 0&0&0&0&0&1&1&-1\end {array} \right] ,\]
\[\tilde{\mathbf P}^{[1]}_2=\frac12\left[ \begin {array}{cccccccc}
\boxed{0}&\boxed{0}&\boxed{0}&\boxed{0}&\boxed{0}&-1&0&-1\\
\boxed{0}&\boxed{0}&\boxed{0}&\boxed{0}&\boxed{0}&-1&-1&0\\
\boxed{0}&\boxed{0}&\boxed{0}&\boxed{0}&\boxed{0}&-1&-1&-1\\
\boxed{0}&\boxed{0}&\boxed{0}&\boxed{0}&\boxed{0}&0&0&-1\\
\boxed{0}&\boxed{0}&\boxed{0}&\boxed{0}&\boxed{0}&0&-1&0\\
-1&-1&-1&0&0&1&1&1\\
0&-1&-1&-1&0&1&1&2\\
-1&0&-1&0&-1&1&0&1
\end {array} \right]
,\]
\[\tilde{\mathbf P}^{[2]}_2=\frac12\left[ \begin {array}{cccccccc}
\boxed{0}&\boxed{0}&\boxed{0}&\boxed{0}&\boxed{0}&1&0&-1\\
\boxed{0}&\boxed{0}&\boxed{0}&\boxed{0}&\boxed{0}&-1&-1&0\\
\boxed{0}&\boxed{0}&\boxed{0}&\boxed{0}&\boxed{0}&-1&-1&1\\
\boxed{0}&\boxed{0}&\boxed{0}&\boxed{0}&\boxed{0}&0&0&-1\\
\boxed{0}&\boxed{0}&\boxed{0}&\boxed{0}&\boxed{0}&0&1&0\\
1&-1&-1&0&0&1&1&-1\\
0&-1&-1&1&0&1&1&-2\\
-1&0&1&0&-1&-1&0&1
\end {array}
 \right]
.\]
Notice that all the boxed entries vanish.
Indeed, as the expansion of the symmetric orthogonalization is
trivial in itself, the pure $\{1,2,3,4,5\}$-terms must vanish in expansion degrees $r\geq1$.

In case $\Psi$ is arbitrary, we can compute  \eqref{U} using   $\tilde Q=\underline{\mathcal O}^{\mathrm{Sy}}(A)$ as the base Clifford system,
we obtain
$\tilde{\mathbf P}^{[1]}_0=[\boxed{\tilde p^{[1]}}]$, $ \tilde{\mathbf P}^{[2]}_0=[\boxed{\tilde p^{[2]}}]$, and
\[\tilde{\mathbf P}^{[1]}_1=\left[ \begin {array}{cccccccc} \boxed{\tilde p_{{1}}^{[1]}}&\boxed{\tilde p_{{2}}^{[1]}}&\boxed{\tilde p_{{3}}^{[1]}}&\boxed{\tilde p_{{4}}^{[1]}}&\boxed{\tilde p_{{5}}^{[1]}}&\tilde p^{[1]}&
\tilde p^{[1]}&\tilde p^{[1]}\end {array} \right],\]
\[\tilde{\mathbf P}^{[2]}_1=\left[ \begin {array}{cccccccc} \boxed{\tilde p_{{1}}^{[2]}}&\boxed{\tilde p_{{2}}^{[2]}}&\boxed{\tilde p_{{3}}^{[2]}}&\boxed{\tilde p_{{4}}^{[2]}}&\boxed{\tilde p_{{5}}^{[2]}}&\tilde p^{[2]}&
\tilde p^{[2]}&-\tilde p^{[2]}\end {array} \right],\]
\[\tilde{\mathbf P}^{[1]}_2= \left[ \begin {smallmatrix}
\boxed{\tilde p_{{1,1}}^{[1]}}&\boxed{\tilde p_{{1,2}}^{[1]}}&\boxed{\tilde p_{{1,3}}^{[1]}}\\
\boxed{\tilde p_{{2,1}}^{[1]}}&\boxed{\tilde p_{{2,2}}^{[1]}}&\boxed{\tilde p_{{2,3}}^{[1]}}\\
\boxed{\tilde p_{{3,1}}^{[1]}}&\boxed{\tilde p_{{3,2}}^{[1]}}&\boxed{\tilde p_{{3,3}}^{[1]}}\\
\boxed{\tilde p_{{4,1}}^{[1]}}&\boxed{\tilde p_{{4,2}}^{[1]}}&\boxed{\tilde p_{{4,3}}^{[1]}}\\
\boxed{\tilde p_{{5,1}}^{[1]}}&\boxed{\tilde p_{{5,2}}^{[1]}}&\boxed{\tilde p_{{5,3}}^{[1]}}\\&\\
-\frac12\tilde p^{[1]}+\frac12 \tilde p_{{1}}^{[1]}&-\frac12\tilde p^{[1]}+\frac12 \tilde p_{{2}}^{[1]}&-\frac12 \tilde p^{[1]}+\frac12 \tilde p_{{3}}^{[1]}\\&\\
\tilde p_{{1}}^{[1]}-\tilde p_{{4}}^{[1]}&-\frac12 \tilde p^{[1]}+\tilde p_{{2}}^{[1]}-\frac12\tilde p_{{3}}^{[1]}&-\frac12 \tilde p^{[1]}-\frac12 \tilde p_{{2}}^{[1]}+\tilde p_{{3}}^{[1]}\\&\\
-\frac12 \tilde p^{[1]}+\frac12\tilde p_{{3}}^{[1]}&\tilde p_{{5}}^{[1]}&-\frac12\tilde p^{[1]}+\frac12 \tilde p_{{1}}^{[1]}\end {smallmatrix}
\right|\ldots\]
\[\ldots\left| \begin {smallmatrix}
\boxed{\tilde p_{{1,4}}^{[1]}}&\boxed{\tilde p_{{1,5}}^{[1]}}&-\frac12  \tilde p^{[1]}+\frac12 \tilde p_{{1}}^{[1]}&\tilde p_{{5}}^{[1]}&-\frac12 \tilde p^{[1]}+\tilde p_{{1}}^{[1]}-\frac12\tilde p_{{3}}^{[1]}\\
\boxed{\tilde p_{{2,4}}^{[1]}}&\boxed{\tilde p_{{2,5}}^{[1]}}&-\frac12 \tilde p^{[1]}+\frac12\tilde p_{{2}}^{[1]}&-\frac12 \tilde p^{[1]}+\frac12 \tilde p_{{3}}^{[1]}&\tilde p_{{2}}^{[1]}-\tilde p_{{4}}^{[1]}\\
\boxed{\tilde p_{{3,4}}^{[1]}}&\boxed{\tilde p_{{3,5}}^{[1]}}&-\frac12 \tilde p^{[1]}+\frac12 \tilde p_{{3}}^{[1]}&-\frac12 \tilde p^{[1]}+\frac12 \tilde p_{{2}}^{[1]}
&-\frac12 \tilde p^{[1]}-\frac12 \tilde p_{{1}}^{[1]}+\tilde p_{{3}}^{[1]}\\
\boxed{\tilde p_{{4,4}}^{[1]}}&\boxed{\tilde p_{{4,5}}^{[1]}}&0&0& -\frac12\tilde p^{[1]}-\frac12 \tilde p_{{2}}^{[1]}+\tilde p_{{4}}^{[1]}\\
\boxed{\tilde p_{{5,4}}^{[1]}}&\boxed{\tilde p_{{5,5}}^{[1]}}&\tilde p_{{5}}^{[1]}&
-\frac12 \tilde p^{[1]}+\frac12 \tilde p_{{1}}^{[1]}&\tilde p_{{5}}^{[1]}\\&\\
0&\tilde p_{{5}}^{[1]}&\frac12 \tilde p^{[1]}-\frac12 \tilde p_{{3}}^{[1]}&\frac12 \tilde p^{[1]}-\frac12\tilde p_{{2}}^{[1]}&\frac12 \tilde p^{[1]}-\frac12\tilde p_{{1}}^{[1]}\\&\\
-\frac12 \tilde p^{[1]}-\frac12 \tilde p_{{1}}^{[1]}+\tilde p_{{4}}^{[1]}&\tilde p_{{5}}^{[1]}&\frac12 \tilde p^{[1]}-\frac12 \tilde p_{{2}}^{[1]}&
\frac12 \tilde p^{[1]}-\frac12\tilde p_{{3}}^{[1]}&\tilde p^{[1]}-\tilde p_{{4}}^{[1]}\\&\\
0&-\frac12\tilde p^{[1]}+\frac12 \tilde p_{{2}}^{[1]}&\frac12 \tilde p^{[1]}-\frac12 \tilde p_{{1}}^{[1]}&-\tilde p_{{5}}^{[1]}&\frac12 \tilde p^{[1]}-\frac12 \tilde p_{{3}}^{[1]}\end {smallmatrix}
\right],\]
\[\tilde{\mathbf P}^{[2]}_2=\left[ \begin {smallmatrix}
\boxed{\tilde p_{{1,1}}^{[2]}}&\boxed{\tilde p_{{1,2}}^{[2]}}&\boxed{\tilde p_{{1,3}}^{[2]}}\\
\boxed{\tilde p_{{2,1}}^{[2]}}&\boxed{\tilde p_{{2,2}}^{[2]}}&\boxed{\tilde p_{{2,3}}^{[2]}}\\
\boxed{\tilde p_{{3,1}}^{[2]}}&\boxed{\tilde p_{{3,2}}^{[2]}}&\boxed{\tilde p_{{3,3}}^{[2]}}\\
\boxed{\tilde p_{{4,1}}^{[2]}}&\boxed{\tilde p_{{4,2}}^{[2]}}&\boxed{\tilde p_{{4,3}}^{[2]}}\\
\boxed{\tilde p_{{5,1}}^{[2]}}&\boxed{\tilde p_{{5,2}}^{[2]}}&\boxed{\tilde p_{{5,3}}^{[2]}}\\&\\
\frac12 \tilde p^{[2]}+\frac12 \tilde p_{{1}}^{[2]}&-\frac12\tilde p^{[2]}+\frac12 \tilde p_{{2}}^{[2]}&-\frac12 \tilde p^{[2]}+\frac12\tilde p_{{3}}^{[2]}\\ &\\
\tilde p_{{1}}^{[2]}-\tilde p_{{4}}^{[2]}&-\frac12 \tilde p^{[2]}+\tilde p_{{2}}^{[2]}-\frac12 \tilde p_{{3}}^{[2]}&-\frac12 \tilde p^{[2]}-\frac12 \tilde p_{{2}}^{[2]}+\tilde p_{{3}}^{[2]}\\ &\\
-\frac12 \tilde p^{[2]}+\frac12 \tilde p_{{3}}^{[2]}&\tilde p_{{5}}^{[2]}&\frac12 \tilde p^{[2]}+\frac12 \tilde p_{{1}}^{[2]}\end
{smallmatrix} \right|\ldots\]
\[\ldots\left| \begin {smallmatrix}
\boxed{\tilde p_{{1,4}}^{[2]}}&\boxed{\tilde p_{{1,5}}^{[2]}}&
\frac12 \tilde p^{[2]}+\frac12\tilde p_{{1}}^{[2]}&\tilde p_{{5}}^{[2]}& -\frac12 \tilde p^{[2]}-\tilde p_{{1}}^{[2]}-\frac12\tilde p_{{3}}^{[2]}\\
\boxed{\tilde p_{{2,4}}^{[2]}}&\boxed{\tilde p_{{2,5}}^{[2]}}&-\frac12\tilde p^{[2]} +\frac12 \tilde p_{{2}}^{[2]}&-\frac12 \tilde p^{[2]}+\frac12\tilde p_{{3}}^{[2]}& -\tilde p_{{2}}^{[2]}-\tilde p_{{4}}^{[2]}\\
\boxed{\tilde p_{{3,4}}^{[2]}}&\boxed{\tilde p_{{3,5}}^{[2]}}&
-\frac12\tilde p_{{}}^{[2]}+\frac12 \tilde p_{{3}}^{[2]}&-\frac12 \tilde p^{[2]}+\frac12\tilde p_{{2}}^{[2]}&\frac12 \tilde p^{[2]}-\frac12 \tilde p_{{1}}^{[2]}-\tilde p_{{3}}^{[2]}\\ \boxed{\tilde p_{{4,4}}^{[2]}}&\boxed{\tilde p_{{4,5}}^{[2]}}&0&0&-\frac12 \tilde p^{[2]}-\frac12 \tilde p_{{2}}^{[2]}-\tilde p_{{4}}^{[2]}\\
\boxed{\tilde p_{{5,4}}^{[2]}}&\boxed{\tilde p_{{5,5}}^{[2]}}&
\tilde p_{{5}}^{[2]}&\frac12 \tilde p^{[2]}+\frac12\tilde p_{{1}}^{[2]}&-\tilde p_{{5}}^{[2]}\\&\\
0&\tilde p_{{5}}^{[2]}&\frac12 \tilde p^{[2]}-\frac12\tilde p_{{3}}^{[2]}&\frac12 \tilde p^{[2]}-\frac12 \tilde p_{{2}}^{[2]}&-\frac12 \tilde p^{[2]}-\frac12\tilde p_{{1}}^{[2]}\\&\\
\frac12 \tilde p^{[2]}-\frac12\tilde p_{{1}}^{[2]}+\tilde p_{{4}}^{[2]}&\tilde p_{{5}}^{[2]}&\frac12
\tilde p^{[2]}-\frac12 \tilde p_{{2}}^{[2]}&\frac12 \tilde p^{[2]}-\frac12 \tilde p_{{3}}^{[2]}&-\tilde p^{[2]}-\tilde p_{{4}}^{[2]}\\&\\
0&-\frac12 \tilde p^{[2]}+\frac12\tilde p_{{2}}^{[2]}&-\frac12\tilde p^{[2]}-\frac12 \tilde p_{{1}}^{[2]}&-\tilde p_{{5}}^{[2]}&\frac12 \tilde p^{[2]}-\frac12 \tilde p_{{3}}^{[2]}\end
{smallmatrix} \right].\]

Hence, we can eliminate the coefficients which are not in the boxed positions.
One can see that the coefficients with any $\{6,7,8\}$-indices decay to linear combinations of coefficients with only
$\{1,2,3,4,5\}$-indices but of lower number.
\end{example}
\begin{example}
Similarly, in the mixed basis, we find that naturality implies
\[\hat{\mathbf P}^{[0]}_0=[\boxed{\hat p^{[0]}}],\qquad
\hat{\mathbf P}^{[0]}_1=\left[ \begin {array}{cccccccc} \boxed{\hat p_{{1}}^{[0]}}&\boxed{\hat p_{{2}}^{[0]}}&\boxed{\hat p_{{3}}^{[0]}}&\boxed{\hat p_{{4}}^{[0]}}&\boxed{\hat p_{{5}}^{[0]}}&0&0&0\end {array} \right],\]
\[\hat{\mathbf P}^{[1]}_0=[\boxed{\hat p^{[1]}}],\qquad
\hat{\mathbf P}^{[1]}_1=\left[ \begin {array}{cccccccc} \boxed{\hat p_{{1}}^{[1]}}&\boxed{\hat p_{{2}}^{[1]}}&\boxed{\hat p_{{3}}^{[1]}}&\boxed{\hat p_{{4}}^{[1]}}&\boxed{\hat p_{{5}}^{[1]}}&\hat p^{[1]}&
\hat p^{[1]}&\hat p^{[1]}\end {array} \right],\]
\[\hat{\mathbf P}^{[2]}_0=[\boxed{\hat p^{[2]}}],\qquad
\hat{\mathbf P}^{[2]}_1=\left[ \begin {array}{cccccccc} \boxed{\hat p_{{1}}^{[2]}}&\boxed{\hat p_{{2}}^{[2]}}&\boxed{\hat p_{{3}}^{[2]}}&\boxed{\hat p_{{4}}^{[2]}}&\boxed{\hat p_{{5}}^{[2]}}&\hat p^{[2]}&
\hat p^{[2]}&-\hat p^{[2]}\end {array} \right],\]
\[\hat{\mathbf P}^{[12]}_0=[\boxed{\hat p^{[12]}}],\qquad
\hat{\mathbf P}^{[12]}_1=\left[ \begin {array}{cccccccc} \boxed{\hat p_{{1}}^{[12]}}&\boxed{\hat p_{{2}}^{[12]}}&\boxed{\hat p_{{3}}^{[12]}}&\boxed{\hat p_{{4}}^{[12]}}&\boxed{\hat p_{{5}}^{[12]}}&0& 2\hat p^{[12]}&0\end {array} \right].\]
\end{example}
\begin{point} \textbf{Transposition invariance \texttt{(Opp)}.}
For any algebra $\mathfrak A$, there is an opposite algebra $\mathfrak A^\opp$.
Using the notation
\[(B_1,B_2)^\opp=(B_1^\opp,B_2^\opp), \]
transposition invariance for scalar and vectorial operations can be expressed as
\[\Xi_{(Q_1,Q_2)}(A_1,A_2)^\opp=\Xi_{(Q_1^\opp,Q_2^\opp)} (A_1^\opp,A_2^\opp);\]
and for pseudoscalar operations,
\[\Xi_{(Q_1,Q_2)}(A_1,A_2)^\opp=-\Xi_{(Q_1^\opp,Q_2^\opp)}(A_1^\opp,A_2^\opp).\]
This invariance means that the FQ operation does not favor right or left.

Regarding the expansions \eqref{f:scalar}, \eqref{f:vec}, \eqref{f:pseudo}, transposition invariance means
that in the splitting basis
\begin{align*}
f_{0}(r_1,\ldots,r_8)&=f_{0}^\opp(r_1,r_2,-r_3,-r_4,r_5,-r_6,r_7,-r_8),\\
f_{1}(r_1,\ldots,r_8)&=f_{1}^\opp(r_1,r_2,r_3,r_4,r_5,-r_6,-r_7,r_8),\\
f_{2}(r_1,\ldots,r_8)&=f_{2}^\opp(r_1,-r_2,-r_3,r_4,r_5,r_6,r_7,r_8),\\
f_{12}(r_1,\ldots,r_8)&=f_{12}^\opp(r_1,-r_2,r_3,-r_4,r_5,r_6,-r_7,-r_8);
\end{align*}
in the mixed basis
\begin{align*}
\hat f_{0}(\hat{r}_1,\ldots,\hat{r}_8)&=\hat f_{0}^\opp(\hat{r}_1,\hat{r}_2,\hat{r}_3,\hat{r}_4,-\hat{r}_5,-\hat{r}_6,-\hat{r}_7,-\hat{r}_8),\\
\hat f_{1}(\hat{r}_1,\ldots,\hat{r}_8)&=\hat f_{1}^\opp(\hat{r}_2,\hat{r}_1,\hat{r}_3,\hat{r}_4,\hat{r}_5,\hat{r}_6,\hat{r}_8,\hat{r}_7),\\
\hat f_{2}(\hat{r}_1,\ldots,\hat{r}_8)&=\hat f_{2}^\opp(-\hat{r}_2,-\hat{r}_1,\hat{r}_3,\hat{r}_4,\hat{r}_5,\hat{r}_6,-\hat{r}_8,-\hat{r}_7),\\
\hat f_{12}(\hat{r}_1,\ldots,\hat{r}_8)&=\hat f_{12}^\opp(-\hat{r}_1,-\hat{r}_2,\hat{r}_3,\hat{r}_4,-\hat{r}_5,-\hat{r}_6,\hat{r}_7,\hat{r}_8);
\end{align*}
in the circular basis
\begin{align*}
\tilde f_{0}(\tilde{r}_1,\ldots,\tilde{r}_8)&=
\tilde f_{0}^\opp(\tilde{r}_1,\tilde{r}_2,\tilde{r}_3,\tilde{r}_5,\tilde{r}_4,-\tilde{r}_6,-\tilde{r}_7,-\tilde{r}_8),\\
\tilde f_{1}(\tilde{r}_1,\ldots,\tilde{r}_8)&=
\tilde f_{1}^\opp(\tilde{r}_2,\tilde{r}_1,\tilde{r}_3,\tilde{r}_4,\tilde{r}_5,\tilde{r}_6,\tilde{r}_8,\tilde{r}_7),\\
\tilde f_{2}(\tilde{r}_1,\ldots,\tilde{r}_8)&=
\tilde f_{2}^\opp(-\tilde{r}_2,-\tilde{r}_1,\tilde{r}_3,\tilde{r}_4,\tilde{r}_5,\tilde{r}_6,-\tilde{r}_8,-\tilde{r}_7),\\
\tilde f_{12}(\tilde{r}_1,\ldots,\tilde{r}_8)&=
\tilde f_{12}^\opp(-\tilde{r}_1,-\tilde{r}_2,\tilde{r}_3,\tilde{r}_5,\tilde{r}_4,-\tilde{r}_6,\tilde{r}_7,\tilde{r}_8).
\end{align*}
(Here $f^\opp$ means that the order of the products is reversed.)
In fact, if $\Xi$ is an FQ operation, then we can define $\Xi^\opp$ by
\[\Xi^\opp(A_1,A_2)=\pm\Xi(A_1^\opp,A_2^\opp)^\opp\]
(minus sign in the pseudoscalar case).
Now, transposition invariance means $\Xi=\Xi^\opp$; and
the RHS of the equations above inform us about the expansion terms of $\Xi^\opp$.

In particular, in the mixed base, transposition invariance implies
\begin{align*}
\begin{bmatrix} \hat p^{[0]}_1&\hat p^{[0]}_2&\hat p^{[0]}_3&\hat p^{[0]}_4&\hat p^{[0]}_5 \end{bmatrix}
&=\begin{bmatrix} \hat p^{[0]}_1&\hat p^{[0]}_2&\hat p^{[0]}_3&\hat p^{[0]}_4&0 \end{bmatrix},\\
\begin{bmatrix} \hat p^{[1]}_1&\hat p^{[1]}_2&\hat p^{[1]}_3&\hat p^{[1]}_4&\hat p^{[1]}_5 \end{bmatrix}
&=\begin{bmatrix} \hat p^{[1]}_1&\hat p^{[1]}_1&\hat p^{[1]}_3&\hat p^{[1]}_4&\hat p^{[1]}_5 \end{bmatrix},\\
\begin{bmatrix} \hat p^{[2]}_1&\hat p^{[2]}_2&\hat p^{[2]}_3&\hat p^{[2]}_4&\hat p^{[2]}_5 \end{bmatrix}
&=\begin{bmatrix} \hat p^{[2]}_1&-\hat p^{[2]}_1&\hat p^{[2]}_3&\hat p^{[2]}_4&\hat p^{[2]}_5 \end{bmatrix},\\
\begin{bmatrix} \hat p^{[12]}_1&\hat p^{[12]}_2&\hat p^{[12]}_3&\hat p^{[12]}_4&\hat p^{[12]}_5 \end{bmatrix}
&=\begin{bmatrix} 0&0&\hat p^{[12]}_3&\hat p^{[12]}_4&0 \end{bmatrix}.
\end{align*}
\end{point}
\begin{point} \textbf{Symmetry \texttt{($\Sigma$${}_2$)}.}
Using the notation
\[(B_1,B_2)^\leftrightarrow=(B_2 ,B_1), \]
rotation invariance for scalar operations can be expressed as
\[\Xi_{(Q_1,Q_2)}(A_1,A_2)=\Xi_{(Q_2,Q_1)}(A_2,A_1),\]
for vectorial  operations
\[\Xi_{(Q_1,Q_2)}(A_1,A_2)^\leftrightarrow=\Xi_{(Q_2,Q_1)}(A_2,A_1),\]
and for pseudoscalar operations
\[\Xi_{(Q_1,Q_2)}(A_1,A_2)=-\Xi_{(Q_2,Q_1)}(A_2,A_1).\]

Regarding the expansions \eqref{f:scalar}, \eqref{f:vec}, \eqref{f:pseudo}, symmetry means
that in the splitting basis
\[f_{0}(r_1,\ldots,r_8)=f_{0}(r_5,r_7,r_6,r_8,r_1,r_3,r_2,r_4),\]
\[f_{2}(r_1,\ldots,r_8)=f_{1}(r_5,r_7,r_6,r_8,r_1,r_3,r_2,r_4),\]
\[f_{12}(r_1,\ldots,r_8)=f_{12}(r_5,r_7,r_6,r_8,r_1,r_3,r_2,r_4);\]
in the mixed basis
\[\hat f_{0}(\hat{r}_1,\ldots,\hat{r}_8)=\hat f_{0}(-\hat{r}_1,\hat{r}_2,\hat{r}_3,-\hat{r}_4,-\hat{r}_5,\hat{r}_6,\hat{r}_7,-\hat{r}_8),\]
\[\hat f_{2}(\hat{r}_1,\ldots,\hat{r}_8)=\hat f_{1}(-\hat{r}_1,\hat{r}_2,\hat{r}_3,-\hat{r}_4,-\hat{r}_5,\hat{r}_6,\hat{r}_7,-\hat{r}_8),\]
\[\hat f_{12}(\hat{r}_1,\ldots,\hat{r}_8)=\hat f_{12}(-\hat{r}_1,\hat{r}_2,\hat{r}_3,-\hat{r}_4,-\hat{r}_5,\hat{r}_6,\hat{r}_7,-\hat{r}_8);\]
and the same scheme works in the circular basis.
If we use the notation
\[(\pi^{\pm}g)(\hat{r}_1,\ldots,\hat{r}_8)=\tfrac12(g(\hat{r}_1,\ldots,\hat{r}_8)\pm g(-\hat{r}_1,\hat{r}_2,\hat{r}_3,-\hat{r}_4,-\hat{r}_5,\hat{r}_6,\hat{r}_7,-\hat{r}_8)),\]
then symmetry can be written as
\[\pi^-\hat f_{0}=0,\qquad
\pi^+\hat f_{2}=\pi^+\hat f_{1},\qquad
\pi^-\hat f_{2}=-\pi^-\hat f_{1},\qquad
\pi^-\hat f_{12}=0,\]
respectively; and the same scheme works in the circular basis.

In particular, in the mixed base, symmetry implies
\begin{align*}
\begin{bmatrix} \hat p^{[0]}_1&\hat p^{[0]}_2&\hat p^{[0]}_3&\hat p^{[0]}_4&\hat p^{[0]}_5 \end{bmatrix}
&=\begin{bmatrix} 0&\hat p^{[0]}_2&\hat p^{[0]}_3&0&0 \end{bmatrix},\\
\begin{bmatrix} \hat p^{[2]}_1&\hat p^{[2]}_2&\hat p^{[2]}_3&\hat p^{[2]}_4&\hat p^{[2]}_5 \end{bmatrix}
&=\begin{bmatrix} -\hat p^{[1]}_1&\hat p^{[1]}_2&\hat p^{[1]}_3&-\hat p^{[1]}_4&-\hat p^{[1]}_5 \end{bmatrix},\\
\begin{bmatrix} \hat p^{[12]}_1&\hat p^{[12]}_2&\hat p^{[12]}_3&\hat p^{[12]}_4&\hat p^{[12]}_5 \end{bmatrix}
&=\begin{bmatrix} 0& \hat p^{[12]}_2&\hat p^{[12]}_3&0&0 \end{bmatrix}.
\end{align*}
\end{point}

\begin{point} \textbf{Orthogonal invariance \texttt{(O${}_2$)}.}
Using the notation
\[\Rot_\phi(B_1,B_2)=(B_1\cos\phi +B_2\sin\phi ,B_2\cos\phi -B_1\sin\phi ), \]
this invariance for (pseudo)scalar operations can be expressed as
\[\Xi_{(Q_1,Q_2)}(A_1,A_2)=\Xi_{\Rot_\phi(Q_1,Q_2)}( \Rot_\phi(A_1,A_2)),\]
and for vectorial  operations
\[\Rot_\phi(\Xi_{(Q_1,Q_2)}(A_1,A_2))=\Xi_{\Rot_\phi(Q_1,Q_2)}( \Rot_\phi(A_1,A_2)).\]
(Strictly speaking, this is rotational invariance, but due to \texttt{(sL)} / \texttt{(vL)} / \texttt{(psL)},
this is equivalent to orthogonal invariance, and, in particular, it is stronger than symmetry.)

In order to describe this invariance, it is better pass to a formal commutative variable $t$ with $t^2=0$, instead of $\phi$.
Then
\[\Rot_t(B_1,B_2)=(B_1 +B_2t ,B_2 -B_1t ). \]
Coefficients of $t$, in calculations like above, correspond to derivations, or vector fields.
According to this, using $t$  is equivalent to using $\phi$. One can also see that considering several
successive such rotations by $t_1,t_2,\ldots$ is equivalent a rotation by $t=t_1+t_2+\ldots$, where the restriction $t^2=0$ is dropped.

Orthogonal invariance is best to be described using the circular basis.
The action $\Rot_t$ induces a derivation $\Delta$  on $\mathfrak F_2$ given by
\[\Delta(Q_1)=Q_2=-Q_1\cdot Q_1Q_2,\quad\Delta(Q_2)=-Q_1=-Q_2 \cdot Q_1Q_2,\]
\[\Delta(\tilde r_1)=-2\tilde r_1Q_1Q_2,\quad\Delta(\tilde r_2)=0,\quad\Delta(\tilde r_3)=0,\quad\Delta(\tilde r_4)=-2\tilde r_4Q_1Q_2,\]
\[\Delta(\tilde r_5)=2\tilde r_5Q_1Q_2,\quad\Delta(\tilde r_6)=0,\quad\Delta(\tilde r_7)=0,\quad\Delta(\tilde r_8)=-2\tilde r_8Q_1Q_2.\]
Then, in terms of the expansions \eqref{f:scalar}, \eqref{f:vec}, \eqref{f:pseudo},
orthogonal invariance is equivalent to
\[\Delta(\tilde f_0(\tilde r_1,\ldots,\tilde r_8))=0,\]
\[\Delta(\tilde f_1(\tilde r_1,\ldots,\tilde r_8)Q_1)=\tilde f_2(\tilde r_1,\ldots,\tilde r_8)Q_2,\]
\[\Delta(\tilde f_2(\tilde r_1,\ldots,\tilde r_8)Q_2)=-\tilde f_1(\tilde r_1,\ldots,\tilde r_8)Q_1,\]
\[\Delta(\tilde f_{12}(\tilde r_1,\ldots,\tilde r_8)Q_1Q_2)=0.\]
Taking into account that $\Delta$ commutes with $\pi^\pm$, and using the notation $\Delta_0=-\tfrac\Delta2$,
\[\Delta_0(\tilde r_1)=\tilde r_1Q_1Q_2,\quad\Delta_0(\tilde r_2)=0,\quad\Delta_0(\tilde r_3)=0,\quad\Delta_0(\tilde r_4)=\tilde r_4Q_1Q_2,\]
\[\Delta_0(\tilde r_5)=-\tilde r_5Q_1Q_2,\quad\Delta_0(\tilde r_6)=0,\quad\Delta_0(\tilde r_7)=0,\quad\Delta_0(\tilde r_8)=\tilde r_8Q_1Q_2;\]
we see that orthogonal invariance means
\[\Delta_0(\tilde f_0(\tilde r_1,\ldots,\tilde r_8))=0,\]
\[\Delta_0(\tilde f_1(\tilde r_1,\ldots,\tilde r_8))(Q_1Q_2)^{-1}=\pi^-\tilde f_1(\tilde r_1,\ldots,\tilde r_8),\]
\[\tilde f_2(\tilde r_1,\ldots,\tilde r_8   )=
\tilde f_1(-\tilde{r}_1,\tilde{r}_2,\tilde{r}_3,-\tilde{r}_4,-\tilde{r}_5,\tilde{r}_6,\tilde{r}_7,-\tilde{r}_8),\]
\[\Delta_0(\tilde f_{12}(\tilde r_1,\ldots,\tilde r_8))=0.\]
This  goes beyond symmetry by
\[\Delta_0(\pi^+\tilde f_0(\tilde r_1,\ldots,\tilde r_8))=0,\]
\[\Delta_0(\pi^+\tilde f_1(\tilde r_1,\ldots,\tilde r_8))=0,\]
\[\Delta_0(\pi^-\tilde f_1(\tilde r_1,\ldots,\tilde r_8))(Q_1Q_2)^{-1}-\pi^-\tilde f_1(\tilde r_1,\ldots,\tilde r_8)=0,\]
\[\Delta_0(\pi^+\tilde f_{12}(\tilde r_1,\ldots,\tilde r_8))=0.\]

The actions on left above act monomially on the expressions of $f_s(\tilde r_1,\ldots,\tilde r_8)$.
Indeed this follows from the shape of the $\Delta_0$ and that $Q_1Q_2$ (anti)commutes with the $\tilde r_i$.
Consequently, invariance means that certain monomial coefficients in $\tilde f_s(\tilde r_1,\ldots,\tilde r_8)$ vanish.
For example, in terms of pseudoscalar invariance, take a monomial $M$ in $\tilde r_1,\ldots,\tilde r_8$.
Then $\Delta_0(M)=\lambda M$. Now, if $\lambda\neq0$ then it means that the coefficient of $M$ in $\tilde f_{12}$ must vanish;
if $\lambda=0$ then there is no restriction.
Hence, one can check the invariance conditions for monomials (in the circular basis) very fast; nevertheless there  are nontrivial patterns.
For example:
\begin{lemma}\label{lem:cir45}
Consider the coefficients $\tilde p_{i_1,\ldots,i_r}^{[s]}$ where $\{i_1,\ldots,i_r\}\subset\{4,5\}$.
Then the cases when rotation invariance does not imply their vanishing are when
\[\Mult^4_{(i_1,\ldots,i_r)}-\Mult^5_{(i_1,\ldots,i_r)}=0\quad\text{and}\quad [s]\text{ is arbitrary}\]
or
\[\Mult^4_{(i_1,\ldots,i_r)}-\Mult^5_{(i_1,\ldots,i_r)}=1\quad\text{and}\quad [s]\in\{[1],[2]\}. \]
(Here $\Mult^x_\iota$ denotes the multiplicity of $x$ in the sequence $\iota$.)
\begin{proof}
One can show that, with this restricted choice of indices, for $\tilde r_\iota= \tilde r_{i_1}\cdots  \tilde r_{i_r}$, the identity
$\Delta_0 \tilde r_\iota= \left(\Mult^4_\iota-\Mult^5_\iota \right)\tilde r_\iota $
holds.
(This is easy to establish if $\iota$ is composed purely from 4's or 5's; then we can show that
inserting 45's or 54's does not change the eigenvalue.)
 From this, and the previous discussion, one can deduce the statement.
\end{proof}
\end{lemma}
\begin{lemma}
For a word $w$ composed from $\{1,2\}$, let $\red w$ denote its shortest reduction by the rules $11=\lambda$  and
$22=\lambda$ (where $\lambda$ is the empty word).

Consider the coefficients $\tilde p_{i_1,\ldots,i_r}^{[s]}$ where $\{i_1,\ldots,i_r\}\subset\{1,2\}$.
Then the cases when rotation invariance does not imply their vanishing are when
\[\red\,(i_1,\ldots,i_r) \in\{\lambda,(2)\}\quad\text{and}\quad [s]\text{ is arbitrary}\]
or
\[ \red\,(i_1,\ldots,i_r)  \in \{(1),(2,1)\}\quad\text{and}\quad [s]\in\{[1],[2]\}. \]
\begin{proof}
With this restricted choice of indices, consider $\tilde r_\iota= \tilde r_{i_1}\cdots  \tilde r_{i_r}$.
One can show that if $\iota$ is reduced (i.~e.~it contains $1$ and $2$ alternating), then
$\Delta_0 \tilde r_\iota= -(-1)^{i_r}\Mult^1_\iota \cdot\tilde r_\iota $ holds (vanishes for $\iota=\lambda$).
Furthermore, one can show that inserting $11$'s or $22$'s does not change the eigenvalue.
From this, and the previous discussion, one can deduce the statement.
\end{proof}
\end{lemma}

In the mixed basis, the differential action is still simple enough,
\[\Delta_0(\hat r_1)=\hat r_1Q_1Q_2,\quad\Delta_0(\hat r_2)=0,\quad\Delta_0(\hat r_3)=0,\quad\Delta_0(\hat r_4)=\hat r_5Q_1Q_2,\]
\[\Delta_0(\hat r_5)=\hat r_4Q_1Q_2,\quad\Delta_0(\hat r_6)=0,\quad\Delta_0(\hat r_7)=0,\quad\Delta_0(\hat r_8)=\hat r_8Q_1Q_2;\]
but it is no longer diagonalized by the monomials, hence the result is more complicated.
In particular, in the mixed base, orthogonal invariance implies
\begin{align*}
\begin{bmatrix} \hat p^{[0]}_1&\hat p^{[0]}_2&\hat p^{[0]}_3&\hat p^{[0]}_4&\hat p^{[0]}_5 \end{bmatrix}
&=\begin{bmatrix} 0&\hat p^{[0]}_2&\hat p^{[0]}_3&0&0 \end{bmatrix},\\
\begin{bmatrix} \hat p^{[1]}_1&\hat p^{[1]}_2&\hat p^{[1]}_3&\hat p^{[1]}_4&\hat p^{[1]}_5 \end{bmatrix}
&=\begin{bmatrix} \hat p^{[1]}_1&\hat p^{[1]}_2&\hat p^{[1]}_3&\hat p^{[1]}_4&\hat p^{[1]}_4 \end{bmatrix},\\
\begin{bmatrix} \hat p^{[2]}_1&\hat p^{[2]}_2&\hat p^{[2]}_3&\hat p^{[2]}_4&\hat p^{[2]}_5 \end{bmatrix}
&=\begin{bmatrix} -\hat p^{[1]}_1&\hat p^{[1]}_2&\hat p^{[1]}_3&-\hat p^{[1]}_4&-\hat p^{[1]}_4 \end{bmatrix},\\
\begin{bmatrix} \hat p^{[12]}_1&\hat p^{[12]}_2&\hat p^{[12]}_3&\hat p^{[12]}_4&\hat p^{[12]}_5 \end{bmatrix}
&=\begin{bmatrix} 0& \hat p^{[12]}_2&\hat p^{[12]}_3&0&0 \end{bmatrix}.
\end{align*}
\end{point}
\begin{point}
Looking for nice FQ operations, the invariance properties above are most desirable,
although they may be invalid for some auxiliary constructions.
Some other, stronger, conditions abstract the basic properties of floating-analytic expansions:

\texttt{(\underline{Biv}')} $\Xi$ is bivariant if for $\theta_1,\theta_2 \approx 1$, the identity
\[\Xi(\theta_1A_1\theta_2,\theta_1A_2\theta_2)=\theta_1\cdot\Xi(A_1,A_2)\cdot\theta_2\]
holds. We talk about symmetric bivariance if this holds with the choice $\theta_1=\theta_2 $.

\texttt{(\underline{Biv}'')} $\Xi$ is antivariant, if for $\theta_1,\theta_2 \approx 1$, the identity
\[\Xi(\theta_1A_1\theta_2,\theta_1A_2\theta_2)=\theta_2^{-1}\cdot\Xi(A_1,A_2)\cdot\theta_1^{-1}\]
holds. We talk about symmetric antivariance if this holds with the choice $\theta_1=\theta_2 $.

\texttt{(\underline{Liv})}  $\Xi$ is left-variant  if for $\theta_1,\theta_2 \approx 1$, the identity
\[\Xi(\theta_1A_1\theta_2,\theta_1A_2\theta_2)=\theta_1\cdot\Xi(A_1,A_2)\cdot\theta_1^{-1}\]
holds. We talk about symmetric left-variance if this holds with the choice $\theta_1=\theta_2 $.

\texttt{(\underline{Riv})}  $\Xi$ is right-variant if for $\theta_1,\theta_2 \approx 1$, the identity
\[\Xi(\theta_1A_1\theta_2,\theta_1A_2\theta_2)=\theta_2^{-1}\cdot\Xi(A_1,A_2)\cdot\theta_2\]
holds. We talk about symmetric right-variance if this holds with the choice $\theta_1=\theta_2 $.
\end{point}

\section{Scaling invariance properties}
\begin{point} \textbf{$\alpha$-Homogeneity \texttt{(\underline{H}${}^\alpha$)}.}
We formulate this invariance infinitesimally. Here $t$ is a commutative formal variable with $t^2=0$.
Then, $\alpha$-homogeneity  can be expressed as
\[(1+\alpha t)\cdot\Xi_{(Q_1,Q_2)}(A_1,A_2)=\Xi_{(Q_1,Q_2)}( (1+t)A_1,(1+ t)A_2)\]
(in the vectorial case the product is taken componentwise).

Here we will use the mixed basis. Then there is a derivation $\Delta_3$ on $(R/Q)$ given by
\[\Delta_3(\hat r_1)=\hat r_1,\quad\Delta_3(\hat r_2)=\hat r_2,\quad\Delta_3(\hat r_3)=\hat r_3+1,\quad\Delta_3(\hat r_4)=\hat r_4,\]
\[\Delta_3(\hat r_5)=\hat r_5,\quad\Delta_3(\hat r_6)=\hat r_6,\quad\Delta_3(\hat r_7)=\hat r_7,\quad\Delta_3(\hat r_8)=\hat r_8\]
(which does not respect the filtration, but anyway).

Then, $\alpha$-homogeneity corresponds to the collection of  relations
\[\hat p^{[s]}_{3,k_1,k_2,\ldots,k_r}+\hat p^{[s]}_{k_1,3,k_2,\ldots,k_r}+\ldots+\hat p^{[s]}_{k_1,k_2,\ldots,k_r,3}
= (\alpha-r)\hat p^{[s]}_{k_1,k_2,\ldots,k_r}.\]
Ultimately, this means that all the $\hat p^{[s]}_{3,k_1,k_2,\ldots,k_r}$ (coefficients starting with index $3$) can be eliminated.
In particular, in the Clifford conservative case, this implies $\hat p_3^{[s]}=\alpha$.
\end{point}
\begin{point} \textbf{$\alpha$-Equiaffinity \texttt{(\underline{E}${}^\alpha$)}.}
Again, we formulate this condition infinitesimally. Let  $t$ be a commutative formal variable with $t^2=0$, and let
\[\mathrm E^\alpha(B_1,B_2):=((1+\alpha t)B_1,(1-\alpha t)B_2).\]

Then, in the vectorial case, $\alpha$-equiaffinity  can be expressed as
\[\mathrm E^\alpha\Xi_{(Q_1,Q_2)}(A_1,A_2)=\Xi_{(Q_1,Q_2)}( (1+t)A_1,(1- t)A_2).\]
In the (pseudo)scalar case, it is defined by
\[(1+\alpha t)\cdot \Xi_{(Q_1,Q_2)}(A_1,A_2)=\Xi_{(Q_1,Q_2)}( (1+t)A_1,(1- t)A_2).\]

Here, $\mathrm E^\alpha$ induces a derivation $\Delta_4$ on $(R/Q)$ given by
\[\Delta_4(\hat r_1)=\hat r_2,\quad\Delta_4(\hat r_2)=\hat r_1,\quad\Delta_4(\hat r_3)=\hat r_4,\quad\Delta_4(\hat r_4)=\hat r_3+1,\]
\[\Delta_4(\hat r_5)=\hat r_6,\quad\Delta_4(\hat r_6)=\hat r_5,\quad\Delta_4(\hat r_7)=\hat r_8,\quad\Delta_4(\hat r_8)=\hat r_7.\]

Then, $\alpha$-equiaffinity corresponds to the collection of  relations
\[\hat p^{[s]}_{4,k_1,k_2,\ldots,k_r}+\hat p^{[s]}_{k_1,4,k_2,\ldots,k_r}+\ldots+\hat p^{[s]}_{k_1,k_2,\ldots,k_r,4}=
(-1^{[s]} )\alpha\hat p^{[s]}_{k_1,k_2,\ldots,k_r}- (\widetilde{\Delta}_4\hat p)^{[s]}_{k_1,k_2,\ldots,k_r};\tag{S4}\label{S4}\]
where $(\widetilde{\Delta}_4\hat p)^{[s]}_{k_1,k_2,\ldots,k_r}$ is an
appropriate linear combination of various $\hat p^{[s]}_{l_1,l_2,\ldots,l_r}$.
In particular, in the Clifford conservative case, this implies $\hat p_4^{[s]}=(-1^{[s]} )\alpha$.
\end{point}
\begin{point} \textbf{$\alpha$-Skew-equiaffinity \texttt{(\underline{CE}${}^\alpha$)}.} \label{po:F}
Let  $t$ be a formal variable with $t^2=0$, but which anticommutes with the $Q_i$ and $R_i$ (hence commutes with the $\hat r_i$) and let
\[\CE^\alpha(B_1,B_2):=((1+\alpha t)B_1,(1-\alpha t)B_2).\]
Then, in the vectorial case, $\alpha$-skew-equiaffinity  can be expressed as
\[\CE^\alpha\Xi_{(Q_1,Q_2)}(A_1,A_2)=\Xi_{(Q_1,Q_2)}( (1+t)A_1,(1- t)A_2).\]
In the (pseudo)scalar case, it is defined by
\[(1+\alpha t)\cdot \Xi_{(Q_1,Q_2)}(A_1,A_2)=\Xi_{(Q_1,Q_2)}( (1+t)A_1,(1- t)A_2).\]

Here, $\CE^\alpha$ induces a derivation  $\Delta_5$  on $(R/Q)$ given by
\[\Delta_5(\hat r_1)=\hat r_7,\quad\Delta_5(\hat r_2)=\hat r_8,\quad\Delta_5(\hat r_3)=\hat r_5,\quad\Delta_5(\hat r_4)=\hat r_6,\]
\[\Delta_5(\hat r_5)=\hat r_3+1,\quad\Delta_5(\hat r_6)=\hat r_4,\quad\Delta_5(\hat r_7)=\hat r_1,\quad\Delta_5(\hat r_8)=\hat r_2.\]

Then, $\alpha$-skew-equiaffinity corresponds to the collection of  relations
\[\hat p^{[s]}_{5,k_1,k_2,\ldots,k_r}+\hat p^{[s]}_{k_1,5,k_2,\ldots,k_r}+\ldots+\hat p^{[s]}_{k_1,k_2,\ldots,k_r,5}=
(-1^{[s]} )\alpha\hat p^{[s]}_{k_1,k_2,\ldots,k_r}- (\widetilde{\Delta}_5\hat p)^{[s]}_{k_1,k_2,\ldots,k_r};\tag{S5}\label{S5}\]
where $(\widetilde{\Delta}_5\hat p)^{[s]}_{k_1,k_2,\ldots,k_r}$ is an
appropriate linear combination of various $\hat p^{[s]}_{l_1,l_2,\ldots,l_r}$.
In particular, in the Clifford conservative case, this implies $\hat p_5^{[s]}=(-1^{[s]} )\alpha$.
\end{point}
\begin{point} \textbf{$\alpha$-Superhomogeneity \texttt{(\underline{SH}${}^\alpha$)}
and  $\alpha$-super-equiaffinity \texttt{(\underline{SE}${}^\alpha$)}.} \label{po:S}
Let $F_1$ and $F_2$ be involutions ($F_i^2=1$) of order $0$, such that $F_i$ anticommutes with $Q_i,R_i$ and commutes with $Q_{1-i},R_{1-i} $, and $F_1$ commutes with $F_2$; and let $t$ be an infinitesimal variable (i. e. of order $1$). This describes an extension of $\mathfrak F_2$.
Then   $t_1=F_1t, t_2=F_2t$ can be considered as formal variables with $t_1^2=t_2^2=t_1t_2=t_2t_1=0$,
such that $t_i$ anticommutes with the $Q_i$ and $R_i$
and commutes with $Q_{1-i}$ and $R_{1-i}$; hence $t_i$ commutes with the $\hat r_j$.
Let
\[\mathrm S^\alpha(B_1,B_2):=((1+\alpha t_1)B_1,(1+\alpha t_2)B_2).\]
$\mathrm S^\alpha$ induces  derivations  on $(R/Q)$: $\frac12(\Delta_1+\Delta_2)$ from $t_1$,  and $\frac12(\Delta_2-\Delta_1)$ from $t_2$,
where
\[\Delta_2(\hat r_1)=\hat r_4,\quad\Delta_2(\hat r_2)=\hat r_3+1,\quad\Delta_2(\hat r_3)=\hat r_2,\quad\Delta_2(\hat r_4)=\hat r_1,\]
\[\Delta_2(\hat r_5)=\hat r_8,\quad\Delta_2(\hat r_6)=\hat r_7,\quad\Delta_2(\hat r_7)=\hat r_6,\quad\Delta_2(\hat r_8)=\hat r_5;\]
and
\[\Delta_1(\hat r_1)=\hat r_3+1,\quad\Delta_1(\hat r_2)=\hat r_4,\quad\Delta_1(\hat r_3)=\hat r_1,\quad\Delta_1(\hat r_4)=\hat r_2,\]
\[\Delta_1(\hat r_5)=\hat r_7,\quad\Delta_1(\hat r_6)=\hat r_8,\quad\Delta_1(\hat r_7)=\hat r_5,\quad\Delta_1(\hat r_8)=\hat r_6.\]

In the vectorial case, $\alpha$-superhomogeneity  is expressed as
\[\frac{1+F_1F_2}{2}\cdot \mathrm S^\alpha\Xi_{(Q_1,Q_2)}(A_1,A_2)=\frac{1+F_1F_2}{2}\cdot\Xi_{(Q_1,Q_2)}( (1+t_1)A_1,(1+t_2)A_2);\]
and, in the (pseudo)scalar case, it can be defined as
\[\frac{1+F_1F_2}{2}\cdot (1+\tfrac\alpha 2 t_1+\tfrac\alpha 2 t_2)\cdot \Xi_{(Q_1,Q_2)}(A_1,A_2)=\Xi_{(Q_1,Q_2)}( (1+t_1)A_1,(1+ t_2)A_2).\]

Then, $\alpha$-superhomogeneity corresponds to the collection of  relations
\[\hat p^{[s]}_{2,k_1,k_2,\ldots,k_r}+\hat p^{[s]}_{k_1,2,k_2,\ldots,k_r}+\ldots+\hat p^{[s]}_{k_1,k_2,\ldots,k_r,2}=
\alpha\hat p^{[s]}_{k_1,k_2,\ldots,k_r}- (\widetilde{\Delta}_2\hat p)^{[s]}_{k_1,k_2,\ldots,k_r};\]
where $(\widetilde{\Delta}_2\hat p)^{[s]}_{k_1,k_2,\ldots,k_r}$ is an
appropriate  linear combination of various $\hat p^{[s]}_{l_1,l_2,\ldots,l_r}$.
In particular, in the Clifford conservative case, this implies $\hat p_2^{[s]}=\alpha$.

In the vectorial case, $\alpha$-super-equiaffinity  is expressed as
\[\frac{1-F_1F_2}{2}\cdot \mathrm S^\alpha\Xi_{(Q_1,Q_2)}(A_1,A_2)=\frac{1-F_1F_2}{2}\cdot\Xi_{(Q_1,Q_2)}( (1+t_1)A_1,(1+t_2)A_2);\]
and, in the (pseudo)scalar case, it can be defined as
\[\frac{1-F_1F_2}{2}\cdot (1+\tfrac\alpha 2 t_1-\tfrac\alpha 2 t_2)\cdot \Xi_{(Q_1,Q_2)}(A_1,A_2)=\Xi_{(Q_1,Q_2)}( (1+t_1)A_1,(1+ t_2)A_2).\]

Then, $\alpha$-super-equiaffinity corresponds to the collection of  relations
\[\hat p^{[s]}_{1,k_1,k_2,\ldots,k_r}+\hat p^{[s]}_{k_1,1,k_2,\ldots,k_r}+\ldots+\hat p^{[s]}_{k_1,k_2,\ldots,k_r,1}=
(-1^{[s]} )\alpha\hat p^{[s]}_{k_1,k_2,\ldots,k_r}- (\widetilde{\Delta}_1\hat p)^{[s]}_{k_1,k_2,\ldots,k_r};\]
where $(\widetilde{\Delta}_1\hat p)^{[s]}_{k_1,k_2,\ldots,k_r}$ is an
appropriate  linear combination of various $\hat p^{[s]}_{l_1,l_2,\ldots,l_r}$.
In particular, in the Clifford conservative case, this implies $\hat p_1^{[s]}=(-1^{[s]} )\alpha$.
\end{point}
At this point, one may wonder if there are similar scaling properties related to certain derivations
$\Delta_6$,  $\Delta_7$, $\Delta_8$.
This is indeed the case, even if not so interesting.
\begin{point} \textbf{$\alpha$-Conjugation-invariance: skew-homogenous \texttt{(\underline{CH}${}^\alpha$)} aspects.}
Let $t$ be an infinitesimal variable anticommuting with the $Q_i,R_i$, as in \ref{po:F}.
Then there is a conjugation action  given by
\[\CH^{\alpha}(X):=(1+\tfrac\alpha2t)X(1-\tfrac\alpha2t).\]

Now, skew-homogeneous $\alpha$-conjugation invariance for $\Xi$ means
\[ \CH^{\alpha}\Xi(A_1,A_2)=\Xi(\CH^1A_1,\CH^1A_2) ,\]
where $\CH^{\alpha}$ is meant componentwise in the vectorial case.

It turns out that $\CH^{\alpha}(X)$  induces a derivation $\Delta_6$ given by
\[\Delta_6(\hat r_1)=\hat r_8,\quad\Delta_6(\hat r_2)=\hat r_7,\quad\Delta_6(\hat r_3)=\hat r_6,\quad\Delta_6(\hat r_4)=\hat r_5,\]
\[\Delta_6(\hat r_5)=\hat r_4,\quad\Delta_6(\hat r_6)=\hat r_3+1,\quad\Delta_6(\hat r_7)=\hat r_2,\quad\Delta_6(\hat r_8)=\hat r_1.\]

Then, skew-homogeneous $\alpha$-conjugation-invariance corresponds to the collection of  relations
\[\hat p^{[s]}_{6,k_1,k_2,\ldots,k_r}+\hat p^{[s]}_{k_1,6,k_2,\ldots,k_r}+\ldots+\hat p^{[s]}_{k_1,k_2,\ldots,k_r,6}=
\alpha\hat p^{[s]}_{k_1,k_2,\ldots,k_r}- (\widetilde{\Delta}_6\hat p)^{[s]}_{k_1,k_2,\ldots,k_r}.\]
In particular, in the Clifford conservative case, it yields $\hat p_6^{[s]}=\alpha$.
On the other hand, the naturality / conjugation invariance property always implies this scaling with $\alpha=0$ in the (pseudo)scalar case and with $\alpha=1$ in the vectorial case.
\end{point}
\begin{point} \textbf{$\alpha$-Conjugation-invariance:
superhomogeneous \texttt{(\underline{CSH}${}^\alpha$)} and super-equiaffine \texttt{(\underline{CSE}${}^\alpha$)} aspects.}
Let $t_1$ and $t_2$ be infinitesimal variables as in  \ref{po:S}.
Then there is a conjugation action  given by
\[\CS^{\alpha}(X):=(1+\tfrac\alpha2 (t_1+t_2))X(1-\tfrac\alpha2 (t_1+t_2)).\]

In the vectorial case, superhomogeneous $\alpha$-conjugation-invariance for $\Xi$ means
\[ \frac{1+F_1F_2}2\CS^{\alpha}\Xi(A_1,A_2)=\frac{1+F_1F_2}2\Xi(\CS^1A_1,\CS^1A_2) ;\]
where $\CS^{\alpha}$ is meant componentwise ; and in the (pseudo)scalar case it is defined by
\[ \frac{1+F_1F_2}2  (1+\tfrac\alpha2t_1+\tfrac\alpha2t_2)\CS^{\alpha}\Xi(A_1,A_2)=\frac{1+F_1F_2}2\Xi(\CS^1A_1,\CS^1A_2) .\]
In the vectorial case, super-equiaffine $\alpha$-conjugation-invariance for $\Xi$ means
\[ \frac{1-F_1F_2}2\CS^{\alpha}\Xi(A_1,A_2)=\frac{1-F_1F_2}2\Xi(\CS^1A_1,\CS^1A_2) ;\]
and in the (pseudo)scalar case it is defined by
\[ \frac{1-F_1F_2}2  (1-\tfrac\alpha2t_1+\tfrac\alpha2t_2)\CS^{\alpha}\Xi(A_1,A_2)=\frac{1-F_1F_2}2\Xi(\CS^1A_1,\CS^1A_2) .\]

There are induced derivations given by
\[\Delta_7(\hat r_1)=\hat r_5,\quad\Delta_7(\hat r_2)=\hat r_6,\quad\Delta_7(\hat r_3)=\hat r_7,\quad\Delta_7(\hat r_4)=\hat r_8,\]
\[\Delta_7(\hat r_5)=\hat r_1,\quad\Delta_7(\hat r_6)=\hat r_2,\quad\Delta_7(\hat r_7)=\hat r_3+1,\quad\Delta_7(\hat r_8)=\hat r_4;\]
and
\[\Delta_8(\hat r_1)=\hat r_6,\quad\Delta_8(\hat r_2)=\hat r_5,\quad\Delta_8(\hat r_3)=\hat r_8,\quad\Delta_8(\hat r_4)=\hat r_7,\]
\[\Delta_8(\hat r_5)=\hat r_2,\quad\Delta_8(\hat r_6)=\hat r_1,\quad\Delta_8(\hat r_7)=\hat r_4,\quad\Delta_8(\hat r_8)=\hat r_3+1.\]

Then superhomogeneous $\alpha$-conjugation-invariance corresponds to the collection of  relations
\[\hat p^{[s]}_{7,k_1,k_2,\ldots,k_r}+\hat p^{[s]}_{k_1,7,k_2,\ldots,k_r}+\ldots+\hat p^{[s]}_{k_1,k_2,\ldots,k_r,7}=
\alpha\hat p^{[s]}_{k_1,k_2,\ldots,k_r}- (\widetilde{\Delta}_7\hat p)^{[s]}_{k_1,k_2,\ldots,k_r};\]
and super-equiaffine $\alpha$-conjugation-invariance corresponds to the collection of  relations
\[\hat p^{[s]}_{8,k_1,k_2,\ldots,k_r}+\hat p^{[s]}_{k_1,8,k_2,\ldots,k_r}+\ldots+\hat p^{[s]}_{k_1,k_2,\ldots,k_r,8}=
(-1^{[s]} )\alpha\hat p^{[s]}_{k_1,k_2,\ldots,k_r}- (\widetilde{\Delta}_8\hat p)^{[s]}_{k_1,k_2,\ldots,k_r}.\]

In particular, in the Clifford conservative case,  superhomogeneity implies
$\hat p_7^{[s]}=\alpha$.
However, the  naturality condition always implies this scaling condition with $\alpha=0$ in the scalar case, $\alpha=1$ in the vectorial case,
and  $\alpha=2$ in the pseudoscalar case.
In the super-equiaffine case, $\hat p_8^{[s]}=(-1^{[s]} )\alpha$ should hold.
Naturality always  implies this scaling condition with $\alpha=0$ in the (pseudo)scalar case and $\alpha=1$ in the vectorial case.
\end{point}
\begin{point}\textbf{Summary and general patterns.}
In this section we have studied some scaling invariances adapted to the mixed basis.
Then every element of the mixed basis corresponds to a scaling invariance.
Regarding their behaviour, the variables of the mixed basis can be:
homogeneous ($\{2,3,6,7\}$) or  equiaffine ($\{1,4,5,8\}$);
ordinary ($\{3,4,5,6\}$) or super ($\{1,2,7,8\}$);
straight ($\{1,2 ,3,4\}$) or skew ($\{5,6 ,7,8\}$).
In fact, every $\hat r_i$ corresponds to a character, a $\pm1$ sequence, according to formula
\eqref{eq:charac}. We can define a  group structure $*$ on the ciphers $\{1,\ldots,8\}$ corresponding to the
postionwise multiplication of the characters.
Here `$3$' turns out to be the the identity element (as it corresponds to the sequence containing only 1's).
This ``character group'' $\mathfrak C$ is isomorphic to $(Z_2)^3$ with multiplication table
\[\begin {array}{c|cccccccc}
*&1&2&3&4&5&6&7&8\\
\hline
1&3&4&1&2&7&8&5&6\\
2&4&3&2&1&8&7&6&5\\
3&1&2&3&4&5&6&7&8\\
4&2&1&4&3&6&5&8&7\\
5&7&8&5&6&3&4&1&2\\
6&8&7&6&5&4&3&2&1\\
7&5&6&7&8&1&2&3&4\\
8&6&5&8&7&2&1&4&3
\end {array} .\]

Then every cipher $i$ induces a derivation $\Delta_i$ on $(R/Q)$ defined by
\[\Delta_i(\hat r_j)=\delta_{i,j}1+\hat r_{i*j}.\]

Using this $\Delta_i$, we can effectively write down the $\hat r_i$-scaling-invariance in terms of expansion coefficients as
\[\hat p^{[s]}_{i,k_1,k_2,\ldots,k_r}+\hat p^{[s]}_{k_1,i,k_2,\ldots,k_r}+\ldots+\hat p^{[s]}_{k_1,k_2,\ldots,k_r,i}=
(-1^{[s]} )^{\delta_{i\in \{1,4,5,8\}}}\alpha\hat p^{[s]}_{k_1,k_2,\ldots,k_r}- (\widetilde{\Delta}_i\hat{p})^{[s]}_{k_1,k_2,\ldots,k_r},\]
where
\[(\widetilde{\Delta}_i \hat{p})^{[s]}_{k_1,k_2,\ldots,k_r}=
\hat{p}^{[s]}_{i*k_1,k_2,\ldots,k_r}+\hat{p}^{[s]}_{k_1,i*k_2,\ldots,k_r}+\ldots+\hat{p}^{[s]}_{k_1,k_2,\ldots,i*k_r}.\]

In the Clifford conservative case $\hat p^{[0]}=1$ / $ \hat p^{[1]}=\hat p^{[2]}=1$ / $\hat p^{[12]}=1$, the
degrees of homogeneity can be recovered from $\hat p^{[s]}_i$.
For $i\in\{6,7,8\}$, however, scaling invariance is always provided by naturality; hence not very interesting,
although it retains some interest for pointed expansions.
\end{point}

\begin{point}
\textbf{ $ \alpha$-circular \texttt{(\underline{XE}${}^\alpha$)} and $\alpha$-skew-circular \texttt{(\underline{CXE}${}^\alpha$)}
invariance.}\label{po:scaltilde}
These are invariances associated to the group algebra elements $\tilde 4=\frac12(\hat 4+\hat 5)\in\mathbb R\mathfrak C$
and $\tilde 5=\frac12(\hat 4-\hat 5)\in\mathbb R\mathfrak C$; where the elements of $\mathfrak C$ are denoted by  $\hat 1,\ldots,\hat 8$,
in order to avoid confusion.

Then $ \alpha$-circular ($\tilde r_4$-scaling) invariance \texttt{(\underline{XE}${}^\alpha$)} is expressed as
$\frac12(\eqref{S4}|_{\alpha=\alpha}+\eqref{S5}|_{\alpha=\alpha});$
and  $\alpha$-skew-circular ($\tilde r_5$-scaling) invariance \texttt{(\underline{CXE}${}^\alpha$)} is expressed as
$\frac12(\eqref{S4}|_{\alpha=\alpha}-\eqref{S5}|_{\alpha=-\alpha});$
where \eqref{S4} and \eqref{S5} refer back to earlier equations.
In the Clifford conservative case,
$\tilde r_4$-scaling invariance implies $\tilde p_4^{[s]}=(-1^{[s]} )\alpha$,
and $\tilde r_5$-scaling invariance implies $\tilde p_5^{[s]}=(-1^{[s]} )\alpha$.
\end{point}
Unfortunately, scaling invariances (in themselves) are rather weak properties.

\section{Hyperscaling}

\begin{point}
We say that the FQ operation $\Xi$ satisfies the hyperscaling property
of type $(J,L,\alpha,\beta)$
in  variable  $\hat r_h$, component $[s]$, if in its expansion relative to the mixed base,
the ``decay'' identities
\begin{align}
\hat p^{[s]}_{h}&=(\alpha+\beta)\hat p^{[s]}\notag\\
\hat p^{[s]}_{h,j,\cdots}&=
\alpha \hat p^{[s]}_{j,\cdots}
-\tfrac12J\hat p^{[s]}_{6*h*j,\cdots}+(-\tfrac12-L)\hat p^{[s]}_{h*j,\cdots}\notag\\
\hat p^{[s]}_{\cdots,i,h}&=
\tfrac12J\hat p^{[s]}_{\cdots,i*h*6}+(-\tfrac12+L)\hat p^{[s]}_{\cdots,i*h}
+\beta\hat p^{[s]}_{\cdots,i}\notag\\
\hat p^{[s]}_{\cdots,i,h,j,\cdots}&=
\tfrac12J\hat p^{[s]}_{\cdots,i*h*6,j,\cdots}+(-\tfrac12+L)\hat p^{[s]}_{\cdots,i*h,j,\cdots}
-\tfrac12J\hat p^{[s]}_{\cdots,i,6*h*j,\cdots}+(-\tfrac12-L)\hat p^{[s]}_{\cdots,i,h*j,\cdots}\notag
\end{align}
hold. Summing up the appropriate terms, we see that hyperscaling of type $(J,L,\alpha,\beta)$
in  variable  $\hat r_h$, component $[s]$ implies an $(\alpha+\beta)$-scaling rule in variable $\hat r_h$, component $[s]$,
\[\hat p^{[s]}_{h,k_1,k_2,\ldots,k_r}+\hat p^{[s]}_{k_1,h,k_2,\ldots,k_r}+\ldots+\hat p^{[s]}_{k_1,k_2,\ldots,k_r,h}=
(\alpha+\beta)\hat p^{[s]}_{k_1,k_2,\ldots,k_r}- (\widetilde{\Delta}_h\hat{p})^{[s]}_{k_1,k_2,\ldots,k_r}.\]

Hyperscaling allows us to eliminate coefficients with indices $h$ in the expansion relative to the mixed base.
A related definition is as follows. We say that  FQ operation $\Xi$ satisfies character degeneracy with $\pm 1$,
if in its expansion relative to the mixed base, the identities
\[\hat p^{[s]}_{\cdots,i,\cdots}=\pm\hat p^{[s]}_{\cdots,i*6,\cdots} \]
hold.
\end{point}

\begin{theorem}
(a) An FQ operation satisfies conjugation invariance, i. e. naturality,
if and only if in its expansion relative to the mixed base, it satisfies
\[\text{hyperscaling of type $(0,0,\alpha_6^{[s]},\beta_6^{[s]})$ in $\hat r_6$, component $[s]$,}\tag{C6}\label{C6}\]
\[\text{ hyperscaling of type $(1,0,\alpha_7^{[s]},\beta_7^{[s]})$ in $\hat r_7$, component $[s]$,}\tag{C7}\label{C7}\]
\[\text{ hyperscaling of type $(-1,0,\alpha_8^{[s]},\beta_8^{[s]})$ in $\hat r_8$, component $[s]$;}\tag{C8}\label{C8}\]
where, for scalar operations,
\[(\alpha^{[0]}_6,\beta^{[0]}_6)=(\tfrac12,-\tfrac12),\qquad
(\alpha^{[0]}_7,\beta^{[0]}_7)=(1,-1),\qquad
(\alpha^{[0]}_8,\beta^{[0]}_8)=(0,0);\]
for vectorial operations,
\[(\alpha^{[1]}_6,\beta^{[1]}_6)=(\tfrac12,\tfrac12),\qquad
(\alpha^{[1]}_7,\beta^{[1]}_7)=(1,0),\qquad
(\alpha^{[1]}_8,\beta^{[1]}_8)=(0,1),\]
\[(\alpha^{[2]}_6,\beta^{[2]}_6)=(\tfrac12,\tfrac12),\qquad
(\alpha^{[2]}_7,\beta^{[2]}_7)=(1,0),\qquad
(\alpha^{[2]}_8,\beta^{[2]}_8)=(0,-1);\]
and for pseudoscalar operations
\[(\alpha^{[12]}_6,\beta^{[12]}_6)=(\tfrac12,-\tfrac12),\qquad
(\alpha^{[12]}_7,\beta^{[12]}_7)=(1,1),\qquad
(\alpha^{[12]}_8,\beta^{[12]}_8)=(0,0).\]

(b) The vectorial FQ operation is bivariant,
if and only if it is symmetrically bivariant, if and only if, in addition to \textrm{(C6)--(C9)}, it satisfies
\[\text{character degeneracy with $+1$.}\tag{C45'}\]
(Symmetric) bivariance is inconsistent for (pseudo)scalar FQ operations except for the identically zero operation.

(b') The vectorial FQ operation is antivariant,
if and only if it is symmetrically antivariant, if and only if, in addition to \textrm{(C6)--(C9)}, it satisfies
\[\text{character degeneracy with $-1$.}\tag{C45''}\]
(Symmetric) antivariance is inconsistent for (pseudo)scalar FQ operations except for the identically zero operation.

(c) The (pseudo)scalar FQ operation is left-variant,
if and only if it is symmetrically left-variant, if and only if , in addition to \textrm{(C6)--(C9)}, it satisfies
\[\text{character degeneracy with $+1$.}\tag{C45'}\]
(Symmetric) left-variance is inconsistent for vectorial FQ operations except for the identically zero operation.

(c') The (pseudo)scalar FQ operation is right-variant,
if and only if it is symmetrically right-variant, if and only if , in addition to \textrm{(C6)--(C9)}, it satisfies
\[\text{character degeneracy with $-1$.}\tag{C45''}\]
(Symmetric) right-variance is inconsistent for vectorial FQ operations except for the identically zero operation.
\begin{proof}\enlargethispage{2mm}
We consider $\Xi$ on the perturbation $(A_1,A_2)=(Q_1+R_1,Q_2+R_2)$ of the Clifford system $(Q_1,Q_2)$.

(a) When we check conjugation invariance formally, it is sufficient to check infinitesimally,
i. e. with respect to elements $1+\theta$, where $\theta\mathfrak A\theta=0$.
Even so, we can decompose $\theta$ into (anti)symmetric parts with respect to $Q_1,Q_2$.
Let $C$ denote the $x\mapsto (1+\theta)x(1-\theta)$ conjugation action.

If $\theta$ commutes with $Q_1$ and $Q_2$, then conjugation invariance holds automatically.

If $\theta$  anticommutes with $Q_1$ and $Q_2$, then let $\hat \varrho_1,\ldots,\hat \varrho_8$ be the
mixed base decomposition of $C(A_1,A_2)$ with respect to $(Q_1,Q_2)$.
By direct computation we find that
\[\hat\rho_i=\hat r_i+2\delta_{i,6}\theta+\theta \hat r_{i*6}+\hat r_{i*6}\theta. \]
Then, in terms of the power series expansion, conjugation invariance means
\[\hat f_s(\hat\rho_1,\ldots,\hat\rho_8)Q^{[s]}=(1+\theta)\hat f_s(\hat r_1,\ldots,\hat r_8)Q^{[s]}(1-\theta).\]
Considering the terms which are of order $1$ in $\theta$,  taking into account the (anti)commutation rules of
$\theta$ and  $Q^{[s]}$, and the noncommutativity of the power series; the equality above translates to (C6).

If $\theta$ anticommutes with $Q_1Q_2$, then it can be assumed that $\theta=\theta_0F_1+\theta_0F_2 $,
where $F_i^2=1$, $F_i$ anticommutes with $Q_i$, $F_i$ commutes with $Q_{1-i}$ and the $\hat r_j$, and $\theta_0$ commutes with $Q_j$.
Again, let $\hat \varrho_1,\ldots,\hat \varrho_8$ be the
mixed base decomposition of $C(A_1,A_2)$ with respect to $(Q_1,Q_2)$.
By direct computation we find that
\[\frac{1+F_1F_2}2\hat\rho_i=
\frac{1+F_1F_2}2\left(\hat r_i+2\delta_{i,7}\theta_0F_1+\theta_0F_1 \hat r_{i*7}+\theta_0F_1 \hat r_{i*2}+\hat r_{i*7}\theta_0F_1-\hat r_{i*2}\theta_0F_1\right) \]
and
\[\frac{1-F_1F_2}2\hat\rho_i=
\frac{1-F_1F_2}2\left(\hat r_i+2\delta_{i,8}\theta_0F_1+\theta_0F_1 \hat r_{i*8}
-\theta_0F_1 \hat r_{i*1}+\hat r_{i*7}\theta_0F_1+\hat r_{i*1}\theta_0F_1\right). \]
Then the equalities
\[\frac{1+F_1F_2}2\hat f_s(\hat\rho_1,\ldots,\hat\rho_8)Q^{[s]}
=\frac{1+F_1F_2}2(1+\theta)\hat f_s(\hat r_1,\ldots,\hat r_8)Q^{[s]}(1-\theta)\]
and
\[\frac{1-F_1F_2}2\hat f_s(\hat\rho_1,\ldots,\hat\rho_8)Q^{[s]}
=\frac{1-F_1F_2}2(1+\theta)\hat f_s(\hat r_1,\ldots,\hat r_8)Q^{[s]}(1-\theta)\]
compared with
\[\frac{1+F_1F_2}2\theta= \frac{1+F_1F_2}2\cdot 2 \theta_0F_1 \qquad \text{and}\qquad \frac{1-F_1F_2}2\theta= 0 \]
yield (C7) and (C8).

(b) When we extend to bivariance, it is sufficient to check infinitesimal bivariance,
i. e. with respect to elements $1+\theta$, where $\theta\mathfrak A\theta=0$; and we start checking out
bivariance with respect to the symmetric bivariance action $B$ given by $x\mapsto (1+\theta)x(1+\theta)$.
Again, we can decompose $\theta$ into (anti)symmetric parts with respect to $Q_1,Q_2$.

If $\theta$ anticommutes with $Q_1$ and $Q_2$, then let $\hat \varrho_1,\ldots,\hat \varrho_8$ be the
mixed base decomposition of $B(A_1,A_2)$ with respect to $(Q_1,Q_2)$.
By direct computation we find that
\[\hat\rho_i=\hat r_i+\theta \hat r_{6*i}-\hat r_{i*6}\theta. \]
Then, in terms of the power series expansion, symmetric bivariance means
\[\hat f_s(\hat\rho_1,\ldots,\hat\rho_8)Q^{[s]}=(1+\theta)\hat f_s(\hat r_1,\ldots,\hat r_8)Q^{[s]}(1+\theta).\]
Considering the terms which are of order $1$ in $\theta$,  taking into account the (anti)commutation rules of
$\theta$ and  $Q^{[s]}$, and the noncommutativity of the power series; the equality above
translates to the collection of relations
\[\hat p^{[s]}_{\cdots,i,6*j,\cdots}=\hat p^{[s]}_{\cdots,i*6,j,\cdots} ;\]
and for (pseudo)scalar operations
\[p^{[s]}=0\,,\qquad \hat p^{[s]}_{6*i,\cdots}=\hat p^{[s]}_{i,\cdots}\,,\qquad \hat p^{[s]}_{\cdots,j*6}=-\hat p^{[s]}_{\cdots,j}\, ;\]
and for vectorial operations
\[\hat p^{[s]}_{6*i,\cdots}=\hat p^{[s]}_{i,\cdots}\,,\qquad \hat p^{[s]}_{\cdots,j*6}=\hat p^{[s]}_{\cdots,j}\,.\]
One can see that for (pseudo)scalar operations this implies the vanishing of all coefficients, while
for vectorial operations, it implies (C45').
For the rest of (infinitesimal) symmetric bivariance, we show that it is equivalent to
\[\text{hyperscaling of type $(0,0,\alpha_3^{[s]},\beta_3^{[s]})$ in $\hat r_3$, component $[s]$,}\tag{C3'}\]
\[\text{ hyperscaling of type $(1,0,\alpha_2^{[s]},\beta_2^{[s]})$ in $\hat r_2$, component $[s]$,}\tag{C2'}\]
\[\text{ hyperscaling of type $(-1,0,\alpha_1^{[s]},\beta_1^{[s]})$ in $\hat r_1$, component $[s]$,}\tag{C1'}\]
where
\[(\alpha^{[1]}_3,\beta^{[1]}_3)=(\tfrac12,\tfrac12),\qquad
(\alpha^{[1]}_2,\beta^{[1]}_2)=(1,0),\qquad
(\alpha^{[1]}_1,\beta^{[1]}_1)=(0,1),\]
\[(\alpha^{[2]}_3,\beta^{[2]}_3)=(\tfrac12,\tfrac12),\qquad
(\alpha^{[2]}_2,\beta^{[2]}_2)=(1,0),\qquad
(\alpha^{[2]}_1,\beta^{[2]}_1)=(0,-1).\]

If $\theta$ commutes with $Q_1$ and $Q_2$, then let $\hat \varrho_1,\ldots,\hat \varrho_8$ be the
mixed base decomposition of $B(A_1,A_2)$ with respect to $(Q_1,Q_2)$.
By direct computation, we find that
\[\hat\rho_i=\hat r_i+2\delta_{i,3}\theta+\theta \hat r_{i}+\hat r_{i}\theta. \]
Then, in terms of the power series expansion, symmetric bivariance means
\[\hat f_s(\hat\rho_1,\ldots,\hat\rho_8)Q^{[s]}=(1+\theta)\hat f_s(\hat r_1,\ldots,\hat r_8)Q^{[s]}(1+\theta).\]
Considering the terms which are of order $1$ in $\theta$,  taking into account the (anti)commutation rules of
$\theta$ and  $Q^{[s]}$, and the noncommutativity of the power series; the equality above translates to (C3').

If $\theta$ anticommutes with $Q_1Q_2$, then it can be assumed that $\theta=\theta_0F_1+\theta_0F_2 $,
where $F_i^2=1$, $F_i$ anticommutes with $Q_i$, $F_i$ commutes with $Q_{1-i}$ and the $\hat r_j$, and $\theta_0$ commutes with $Q_j$.
Again, let $\hat \varrho_1,\ldots,\hat \varrho_8$ be the
mixed base decomposition of $B(A_1,A_2)$ with respect to $(Q_1,Q_2)$.
By direct computation we find that
\[\frac{1+F_1F_2}2\hat\rho_i=
\frac{1+F_1F_2}2\left(\hat r_i+2\delta_{i,2}\theta_0F_1+\theta_0F_1 \hat r_{i*7}+\theta_0F_1 \hat r_{i*2}-\hat r_{i*7}\theta_0F_1+\hat r_{i*2}\theta_0F_1\right) \]
and
\[\frac{1-F_1F_2}2\hat\rho_i=
\frac{1-F_1F_2}2\left(\hat r_i+2\delta_{i,8}\theta_0F_2-\theta_0F_2 \hat r_{i*8}
+\theta_0F_2 \hat r_{i*1}+\hat r_{i*7}\theta_0F_2+\hat r_{i*1}\theta_0F_2\right). \]
Then the equalities
\[\frac{1+F_1F_2}2\hat f_s(\hat\rho_1,\ldots,\hat\rho_8)Q^{[s]}=\frac{1+F_1F_2}2(1+\theta)\hat f_s(\hat r_1,\ldots,\hat r_8)Q^{[s]}(1-\theta)\]
and
\[\frac{1-F_1F_2}2\hat f_s(\hat\rho_1,\ldots,\hat\rho_8)Q^{[s]}=\frac{1-F_1F_2}2(1+\theta)\hat f_s(\hat r_1,\ldots,\hat r_8)Q^{[s]}(1-\theta)\]
yield (C2') and (C1').

However, under character degeneracy with $+1$,
(C3') is equivalent to (C6); (C2') is equivalent to  (C7), (C1') is equivalent to  (C8); so
symmetric bivariance also implies conjugation invariance.
(Infinitesimal)  conjugation invariance and (infinitesimal) symmetric bivariance, however, implies full (infinitesimal)  bivariance.
This argument also shows that  (C45') is sufficient to provide bivariance in addition to (C6)--(C8).

(b') We can proceed in similar manner.
If $\theta$ anticommutes with $Q_1$ and $Q_2$, then let $\hat \varrho_1,\ldots,\hat \varrho_8$ be the
mixed base decomposition of $B(A_1,A_2)$ with respect to $(Q_1,Q_2)$.
We can apply the same analysis as before.
Then, in terms of the power series expansion, symmetric antivariance means
\[\hat f_s(\hat\rho_1,\ldots,\hat\rho_8)Q^{[s]}=(1-\theta)\hat f_s(\hat r_1,\ldots,\hat r_8)Q^{[s]}(1-\theta).\]
Considering the terms which are of order $1$ in $\theta$,  taking into account the (anti)commutation rules of
$\theta$ and  $Q^{[s]}$, and the noncommutativity of the power series; the equality above
translates to the collection of relations
\[\hat p^{[s]}_{\cdots,i,6*j,\cdots}=\hat p^{[s]}_{\cdots,i*6,j,\cdots} ;\]
and for (pseudo)scalar operations
\[p^{[s]}=0\,,\qquad \hat p^{[s]}_{6*i,\cdots}=-\hat p^{[s]}_{i,\cdots}\,,\qquad \hat p^{[s]}_{\cdots,j*6}=\hat p^{[s]}_{\cdots,j}\, ;\]
and for vectorial operations
\[\hat p^{[s]}_{6*i,\cdots}=-\hat p^{[s]}_{i,\cdots}\,,\qquad \hat p^{[s]}_{\cdots,j*6}=-\hat p^{[s]}_{\cdots,j}\,.\]
One can see that for (pseudo)scalar operations, this implies the vanishing of all coefficients, while
for vectorial operations, it implies (C45'').
Then the further aspects of symmetric antivariance can be encoded by
\[\text{hyperscaling of type $(0,0,-\alpha_3^{[s]},-\beta_3^{[s]})$ in $\hat r_3$, component $[s]$,}\tag{C3''}\]
\[\text{ hyperscaling of type $(1,0,-\alpha_2^{[s]},-\beta_2^{[s]})$ in $\hat r_2$, component $[s]$,}\tag{C2''}\]
\[\text{ hyperscaling of type $(-1,0,-\alpha_1^{[s]},-\beta_1^{[s]})$ in $\hat r_1$, component $[s]$.}\tag{C1''}\]
The rest of the argument is analogous.

(c) Again, let $\hat \varrho_1,\ldots,\hat \varrho_8$ be the
mixed base decomposition of $B(A_1,A_2)$ with respect to $(Q_1,Q_2)$.
In terms of the power series expansion, symmetric antivariance means
\[\hat f_s(\hat\rho_1,\ldots,\hat\rho_8)Q^{[s]}=(1+\theta)\hat f_s(\hat r_1,\ldots,\hat r_8)Q^{[s]}(1-\theta).\]
Then we can proceed as in (b).

(c') This is analogous to (b').
\end{proof}
\end{theorem}
\begin{remark}
At first sight, the previous theorem is just plainly more informative the argument in \ref{po:conj}, but this is not completely so.
An advantage of \ref{po:conj} is that it proves that using (C6)--(C8) we are led to an unambiguous reduction
in term of the indices $1,2,3,4,5$ (in the mixed and circular bases).
\end{remark}
\begin{cor} One can compute the coefficients of the expansion of $\underline{\mathcal O}^{\mathrm{Sy}}$ in the
mixed basis, recursively, using

(i) $\hat p^{[1]}=\hat p^{[2]}=1$;

(ii) $\hat p^{[s]}_{\iota_1,\ldots,\iota_r}=0$ if $\{\iota_1,\ldots,\iota_r\}$ is a nonempty subset of $\{1,\ldots,5\}$;

(iii) the vectorial hyperscaling rules \eqref{C6},  \eqref{C7}, \eqref{C8} of the previous theorem.
\qed
\end{cor}
\begin{point}
One can prove similar statements regarding the splitting and circular bases, too.
Then, instead of $*$, one should use more complicated incidence matrices of indices.
In particular, in case of the circular basis,  character degeneracy yields
$\tilde p^{[s]}_{\cdots,i,\cdots}=\pm\tilde p^{[s]}_{\cdots,i*6,\cdots} $ for $i\in\{1,2,3,6,7,8\}$;
and, $\tilde p^{[s]}_{\cdots,5,\cdots}=0$ in the $+1$ case, and $\tilde p^{[s]}_{\cdots,4,\cdots}=0$ in the $-1$ case.
\end{point}
\begin{theorem} \label{thm:biv1}
(b) If $\Xi$ is a bivariant vectorial FQ operation, then it is determined by the
coefficients $\tilde p^{[s]}_{4,\ldots,4}$, $s\in\{1,2\}$,  which can be prescribed arbitrarily.

(b')  If $\Xi$ is an antivariant vectorial FQ operation, then it is determined by the
coefficients $\tilde p^{[s]}_{5,\ldots,5}$,  $s\in\{1,2\}$,  which can be prescribed arbitrarily.

(c) If $\Xi$ is a left-variant (pseudo)scalar FQ operation, then it is determined by the
coefficients $\tilde p^{[s]}_{4,\ldots,4}$, ($s=0$ or $s=12$)  which can be prescribed arbitrarily.

(c') If $\Xi$ is a right-variant (pseudo)scalar FQ operation, then it is determined by the
coefficients $\tilde p^{[s]}_{5,\ldots,5}$, ($s=0$ or $s=12$)  which can be prescribed arbitrarily.

(The statement is also true using the mixed base, and then we do not even have worry whether 4 or 5 should be used as indices.)
\begin{proof}
(b) First considering the mixed basis, due to character degeneracy, we can pass to the denegerate base, where
$\hat 1=\hat 8, \hat 2=\hat 7, \hat 3=\hat 7, \hat 4=\hat 5  $.
Then all we have to do is to impose (C6)--(C8).
When we do this we reduce everything to coefficients $\hat p_{4,\ldots,4}^{[s]}$.
The reduction to these terms using  (C6)--(C8) leads to unambiguous results;
or, said differently, it imposes no  relations between the $\hat p_{4,\ldots,4}^{[s]}$, because
we have that much freedom ($r$ linear degrees of freedom up to order $r$)
even in a floating analytic expansion (discussed in \cite{L1} in the bivariant and antivariant cases,
and left to the reader in the left- and right-variant cases).
When we pass to the circular basis, we see that all the coefficients $\tilde p^{[s]}_{\cdots,5,\cdots}$ must vanish, so
the setting is in fact reduced to the coefficients $\tilde p_{4,\ldots,4}^{[s]}$.

(b') First considering the mixed basis, due to character degeneracy, we can pass to the denegerate base,
similarly. The accounting is just a little bit trickier due to the sign changes
$\hat 1=-\hat 8, \hat 2=-\hat 7, \hat 3=-\hat 7, \hat 4=-\hat 5  $.
Again we arrive to reduction to the coefficients $\hat p_{4,\ldots,4}^{[s]}$.
When we pass to the circular basis, we see that all the coefficients $\tilde p^{[s]}_{\cdots,4,\cdots}$ must vanish, so
the setting is, in fact, reduced to the coefficients $\tilde p_{5,\ldots,5}^{[s]}$.

(c) and (c') can be proven analogously.
\end{proof}
\end{theorem}
\begin{example}
In \cite{L1} we have introduced the conform orthogonization procedure $\mathcal O^{\mathrm{fSy}}$.
We define the anticonform orthogonization procedure ${\mathcal O}^{\mathrm{afSy}}$ such that
\[{\mathcal O}^{\mathrm{afSy}}(A_1,A_2)_i:=-(\mathcal O^{\mathrm{fSy}}(A_1,A_2)_i)^{-1};\]
i. e., relative to $\mathcal O^{\mathrm{fSy}}$, we take the inverse times $-1$ in every component.
$\mathcal O^{\mathrm{fSy}}$ and ${\mathcal O}^{\mathrm{afSy}}$ are vectorial FQ operations.

We define the pseudoscalar operation left axis  $\mathcal A_L$ by
\[\mathcal A_L(A_1,A_2):=-\pol A_1A_2^{-1}=\pol A_2A_1^{-1};\]
and the right axis  $\mathcal A_R$ by
\[\mathcal A_R(A_1,A_2):=-\pol A_1^{-1}A_2=\pol A_2^{-1}A_1.\]

These operations above are, in fact, analytic FQ operations. One can show that
$\mathcal O^{\mathrm{fSy}}$ is bivariant,  ${\mathcal O}^{\mathrm{afSy}}$ is antivariant,  $\mathcal A_L$ is left-variant,
$\mathcal A_R$ is right-variant; and they are all orthogonal invariant and Clifford conservative.
In what follows, we will consider their formal restrictions (in notation: underlined).
We remark that in their (formal) expansion
\begin{align}
&\mathcal O^{\mathrm{fSy}}:&&\tilde {\mathbf P}^{[1]}_0=[1]&&\tilde {\mathbf P}^{[1]}_1=\left[ \begin {array}{cccccccc} 1&1&1&0&0&1&1&1\end {array}\right],\notag\\
&&&\tilde {\mathbf P}^{[2]}_0=[1]&&\tilde {\mathbf P}^{[2]}_1=\left[ \begin {array}{cccccccc} -1&1&1&0&0&1&1&-1\end {array}\right];\notag\\
&{\mathcal O}^{\mathrm{afSy}}:&&\tilde {\mathbf P}^{[1]}_0=[1]&&\tilde{\mathbf P}^{[1]}_1=\left[ \begin {array}{cccccccc} -1&-1&-1&0&0&1&1&1\end {array}\right],\notag\\
&&&\tilde {\mathbf P}^{[2]}_0=[1]&&\tilde {\mathbf P}^{[2]}_1=\left[ \begin {array}{cccccccc} 1&-1&-1&0&0&1&1&-1\end {array}\right];\notag\\
&\mathcal A_L:&&\tilde {\mathbf P}^{[12]}_0=[1]&&\tilde {\mathbf P}^{[12]}_1=\begin{bmatrix}0&2&0&0&0&0&2&0\end{bmatrix};\notag\\
&\mathcal A_R:&&\tilde {\mathbf P}^{[12]}_0=[1]&&\tilde {\mathbf P}^{[12]}_1=\begin{bmatrix}0&-2&0&0&0&0&2&0\end{bmatrix}.\notag
\end{align}
\end{example}
\begin{theorem} \label{thm:biv2}
(a) If $\Xi$ is a left-variant or right-variant scalar  FQ operation
which is orthogonal invariant such that $\tilde p^{[0]}=1$, then \[\Xi=\underline{1}.\]

(b) If $\Xi$ is a bivariant vectorial FQ operation which is orthogonal invariant such that $\tilde p^{[1]}=\tilde p^{[2]}=1$, then
\[\Xi=\alpha\cdot\underline{\Id}+(1-\alpha)\cdot\underline{\mathcal O}^{\mathrm{fSy}},\]
where $\tilde p^{[1]}_4=-\tilde p^{[2]}_4=\alpha$.

(b')  If $\Xi$ is an antivariant vectorial FQ operation which is orthogonal invariant such that $\tilde p^{[1]}=\tilde p^{[2]}=1$, then
\[\Xi=\underline{{\mathcal O}}^{\mathrm{afSy}}.\]

(c) If $\Xi$ is a left-variant pseudoscalar FQ operation which is orthogonal invariant such that $\tilde p^{[12]}=1$, then
\[\Xi=\underline{\mathcal A}_L.\]

(c')  If $\Xi$ is a right-variant pseudoscalar FQ operation which is orthogonal invariant such that $\tilde p^{[12]}=1$, then
\[\Xi=\underline{\mathcal A}_R.\]
\begin{proof}
(a) If $\Xi$ is a left-variant scalar operation, then according to the previous theorem,
the operation depends only on the collection of coefficients $\tilde p^{[s]}_{4,\ldots,4}$.
According to our prescriptions, $\tilde p^{[0]}=1$; while, according to Lemma \ref{lem:cir45}, orthogonal invariance implies
that  $\tilde p^{[s]}_{4,\ldots,4}=0$ if the number of lower indices is bigger then zero.
This implies that there is at most one such operation.
However, the choice $\Xi=1$ satisfies the requirements for such an operation.

A very similar argument applies if  $\Xi$ is a right-variant scalar operation, and in cases (b'), (c), (c').
Case (b) is a little bit different, because Lemma \ref{lem:cir45} does not tell about the vanishing of
$\tilde p^{[1]}_{4}=-\tilde p^{[2]}_{4}$ (the equality follows from symmetry).
And indeed, the linear combination given in the statement of (b)
allows arbitrary choice for $\tilde p^{[1]}_{4}=-\tilde p^{[2]}_{4}$.
\end{proof}
\end{theorem}
\begin{point}
Sometimes it is useful to consider hyperscaling with respect to the variables $\tilde r_4$ and $\tilde r_5$.
This means with respect to
$\tilde 4=\frac12(\hat 4+\hat 5)\in\mathbb R\mathfrak C$
and $\tilde 5=\frac12(\hat 4-\hat 5)\in\mathbb R\mathfrak C$;
the corresponding equations are basically the sums and differences of  $\hat r_4$ and $\hat r_5$-scaling equations
similarly to as in \ref{po:scaltilde}.
\end{point}
Hyperscaling constraints exhibit a structured and nontrivial behaviour, which we cannot discuss here.
However, if it is said that scaling constraints are too weak, then it must be said that
hyperscaling constraints are too restrictive.

\section{Floating Clifford conservativity}

\begin{point}
$\texttt{(vC')}$ Floating Clifford conservativity:
A vectorial FQ operation $\Xi$ satisfies this property if it acts trivially on floating Clifford systems.

$\texttt{(vC'')}$ Floating Clifford anticonservativity:
A vectorial FQ operation $\Xi$ satisfies this property if it inverts floating Clifford systems, and multiplies them by $-1$.

In order to deal with floating Clifford (anti)conservitivity, we use the following
\end{point}

\begin{theorem}
Suppose that
$(Q_1,Q_2)$ is the symmetric orthogonalization of $(A_1,A_2)=(Q_1+R_1,Q_2+R_2)$, i. e.
$\tilde{r}_6=\tilde{r}_7=\tilde{r}_8=0$ in the circular base.

(a) We claim: If $ (A_1,A_2)$ is floating Clifford system, then $\tilde r_4$ and $\tilde r_5$ can be expressed
from $\tilde r_1,\tilde r_2,\tilde r_3  $, by some fixed explicit power series $\tilde F_4,\tilde F_5$:
\[\tilde r_4=\tilde F_4(\tilde r_1,\tilde r_2,\tilde r_3 )\equiv\tilde r_2\tilde r_1+O((\tilde r_1,\tilde r_2,\tilde r_3)^{\geq3});\]
\[\tilde r_5=\tilde F_5(\tilde r_1,\tilde r_2,\tilde r_3 )\equiv\tilde r_1\tilde r_2+O((\tilde r_1,\tilde r_2,\tilde r_3)^{\geq3}).\]

(b) Conversely, in the general case, the terms $\tilde r_4$ and $\tilde r_5$ can be replaced by
$\tilde F_4(\tilde r_1,\tilde r_2,\tilde r_3 )$
and  $\tilde F_5(\tilde r_1,\tilde r_2,\tilde r_3 )$ in order to yield a floating Clifford system.

Similar statement holds in  the mixed base.
\begin{proof}
(a) From the equations $(A_1A_2^{-1})^2=-1$ and $(A_2^{-1}A_1)^2=-1$ we obtain
\[ \tilde r_4= \tilde r_2\tilde r_1-\tilde r_2\tilde r_4+\tilde r_4\tilde r_2+\tilde r_4\tilde r_3
-2\tilde r_4\tilde r_4 +O((\tilde r_1,\tilde r_2,\tilde r_3,\tilde r_4,\tilde r_5)^{\geq3})\]
and
\[ \tilde r_5=\tilde r_1\tilde r_2+\tilde r_2\tilde r_5+\tilde r_3\tilde r_5-\tilde r_5\tilde r_2 -2\tilde r_5\tilde r_5
 +O((\tilde r_1,\tilde r_2,\tilde r_3,\tilde r_4,\tilde r_5)^{\geq3}).\]
 Iterating these expressions, we find expressions  $\tilde F_4,\tilde F_5$ as indicated.

(b) Consider the algebra generated by
$\underline{Q}_1,\underline{Q}_2,\underline{\tilde r}_1, \underline{\tilde r}_2,\underline{\tilde r}_3 $
subject to the appropriate commutation relations.
Let $(A_1,A_2)=(1+\underline{\tilde r}_2+\frac12\underline{\tilde r}_3)\cdot( \underline{Q}_1,\underline{Q}_2)\cdot(1+\underline{\tilde r}_2+\frac12\underline{\tilde r}_3)$,
which is a floating Clifford system.
Let $(Q_1,Q_2)$ be its symmetric orthogonalization.
Then one finds that
\[\tilde r_i=\underline{\tilde r}_i+O((\underline{\tilde r}_1,\underline{\tilde r}_2,\underline{\tilde r}_3)^{\geq2})
\qquad \text{for}\quad i\in\{1,2,3\}.\]
This implies that $\underline{\tilde r}_1, \underline{\tilde r}_2,\underline{\tilde r}_3$
and $\tilde r_1, \tilde r_2,\tilde r_3$ can be expressed from each other;
which proves that having a floating Clifford system implies no nontrivial relations for $\tilde r_1, \tilde r_2,\tilde r_3$,
and in, fact, their free prescribability.
\end{proof}
\end{theorem}
The process of (b) yields, in fact, an FQ operation $\underline{\mathcal F}^{\mathrm{Sy}}$ producing floating Clifford systems,
compatible with $\underline{\mathcal O}^{\mathrm{Sy}}$ (hence different from $\underline{\mathcal O}^{\mathrm{fSy}}$,
and therefore not bivariant).
\begin{theorem}
For a conjugation-invariant vectorial FQ operation, floating Clifford (anti-) conservativity is equivalent to a
collections of relations
\[\tilde p^{[s]}_{\iota_1,\ldots,\iota_r}=\text{ an inhomogeneous linear expression
of $\tilde p^{[s]}_{\varkappa_1,\ldots,\varkappa_h}$'s  with $h<r$,}\]
where $\{\iota_1,\ldots,\iota_r\}\subset\{1,2,3\}$, $s\in\{1,2\}$. Similar statement holds in  the mixed base.
\begin{proof} We know that a conjugation-invariant vectorial FQ operation can be encoded by
formal power series $\tilde f_s(\tilde r_1,\ldots, \tilde r_5)$  with respect to the symmetric orthogonalization $(Q_1,Q_2)$.
Considering $(A_1,A_2)=((1+\tilde r_1+\tilde r_2 +\tilde r_3 +\tilde r_4)Q_1,(1-\tilde r_1+\tilde r_2 +\tilde r_3 -\tilde r_4)Q_2)$,
we see that floating Clifford conservativity ($+$) and  anticonservativity ($-$) can be expressed by
\[\tilde f_1(\tilde r_1, \tilde r_2,\tilde r_3,\tilde F_4(\tilde r_1,\tilde r_2,\tilde r_3 ),\tilde F_5(\tilde r_1,\tilde r_2,\tilde r_3 ) )
=\pm( (1+\tilde r_1+\tilde r_2 +\tilde r_3 +\tilde F_4(\tilde r_1,\tilde r_2,\tilde r_3 ))Q_1)^{\pm1}Q_1^{-1},\]
\[\tilde f_2(\tilde r_1, \tilde r_2,\tilde r_3,\tilde F_4(\tilde r_1,\tilde r_2,\tilde r_3 ),\tilde F_5(\tilde r_1,\tilde r_2,\tilde r_3 ) )
=\pm( (1-\tilde r_1+\tilde r_2 +\tilde r_3 -\tilde F_4(\tilde r_1,\tilde r_2,\tilde r_3 ))Q_2)^{\pm1}Q_2^{-1}.\]
As $\tilde F_4,\tilde F_5$ are of higher degree,  the equations expand as in the statement.
\end{proof}
\end{theorem}
We see that floating Clifford (anti)conservativity is a quite weak property ($2\cdot3^r$ constraints compared to $2\cdot5^r$
free parameters) although not trivial.

\section{Clifford productivity}
\begin{point}
$\texttt{(CP)}$ Clifford productivity:
A vectorial FQ operation $\Xi$ satisfies this property if it produces Clifford systems.
A pseudoscalar FQ operation $\Xi$ satisfies this property if it produces skew-involutions.

$\texttt{(CP')}$ Floating Clifford productivity:
A vectorial FQ operation $\Xi$ satisfies this property if it produces floating Clifford systems.

One can deal with Clifford productivity as follows.
First, it is reasonable to restrict to the Clifford conservative case; then one can use the
fact that Clifford systems close to each other are conjugates of each other
(in case of floating Clifford systems: translates of each other).
Done carefully, one can  organize the conjugation scheme such that it provides existence and unicity at the same time.
Such an analysis was already considered in \cite{L1}, here we give a  more detailed account.
\end{point}
\begin{conven}
Suppose that $\Xi$ is an FQ operation. If $\Xi$ is (pseudo)scalar operation, then $\Xi^{-1}$ denotes the inverse with respect
to $1$, regarding multiplication in $\mathfrak A$; if $\Xi$ is vectorial operation, then $\Xi^{-1}$ denotes the inverse with respect
to $\Id$, regarding composition of FQ operations (with respect to the same base point as of $\Xi$).

We will use the notation $\sFQ, \vFQ, \psFQ$ for the spaces of scalar, vectorial, or pseudoscalar FQ operations respectively.
We will use the notation $\sFQ^{\leq r},\sFQ^{\langle r\rangle},\sFQ^{\geq r}$
for those scalar FQ operations whose expansion terms in $(R/Q)$ are in degrees $\leq r$, exactly in degree $r$, or in degrees $\geq r$, respectively;
we use similar notion in the other cases, too.

In the case of a concrete FQ operation $\Xi$, let $\Xi^{\leq r},\Xi^{\langle r\rangle},\Xi^{\geq r}$ denote
those FQ operations which we obtain from $\Xi$ by restricting its expansion to the indicated orders.
\end{conven}
\begin{lemma} Consider the maps
\[\xymatrix{\ar@(ul,dl)_{\eta}}\xymatrix{
\sFQ^{\langle r\rangle}\ar[r]_{\lambda}&\ar@/_/[l]_{\varkappa}
\vFQ^{\langle r\rangle}\ar[r]_{\colambda}&\ar@/_/[l]_{\covarkappa}
\sFQ^{\langle r\rangle}&\hspace{-1cm}\oplus\psFQ^{\langle r\rangle}\oplus\sFQ^{\langle r\rangle}
}\xymatrix{\ar@(dr,ur)_{\coeta}}\]
given by
\begin{align*}
\eta^{\langle r\rangle}:&U^{\langle r\rangle}\mapsto (U^{\langle r\rangle}){}^{0}_{Q_1}{}^{0}_{Q_2};\\
\lambda^{\langle r\rangle}:&U^{\langle r\rangle}\mapsto U^{\langle r\rangle}\cdot(Q_1,Q_2)-(Q_1,Q_2)\cdot U^{\langle r\rangle};\\
\varkappa^{\langle r\rangle}:&(V_1^{\langle r\rangle},V_2^{\langle r\rangle}) \mapsto
\frac12(V_1^{\langle r\rangle}Q_1^{-1})^{10}_{Q}+\frac14(V_1^{\langle r\rangle}Q_1^{-1})^{11}_{Q}
+\frac12(V_2^{\langle r\rangle}Q_2^{-1})^{01}_{Q}+\frac14(V_2^{\langle r\rangle}Q_2^{-1})^{11}_{Q};\\
\colambda^{\langle r\rangle}:&(V_1^{\langle r\rangle},V_2^{\langle r\rangle}) \mapsto
([V_1^{\langle r\rangle},Q_1]_+  ,
[V_2^{\langle r\rangle},Q_1]_++[V_1^{\langle r\rangle},Q_2]_+
,[V_2^{\langle r\rangle},Q_2]_+ );\\
\covarkappa^{\langle r\rangle}:&(W_1^{\langle r\rangle},W_{12}^{\langle r\rangle},W_2^{\langle r\rangle}) \mapsto\\
&\left(-\frac14[W_1^{\langle r\rangle},Q_1]_+-\frac18[W_{12}^{\langle r\rangle},Q_1Q_2]_+\cdot Q_1  ,
-\frac14[W_2^{\langle r\rangle},Q_2]_++\frac18[W_{12}^{\langle r\rangle},Q_1Q_2]_+\cdot Q_2
 \right);\\
\coeta^{\langle r\rangle}:&(W_1^{\langle r\rangle},W_{12}^{\langle r\rangle},W_2^{\langle r\rangle}) \mapsto\\
&\left( (W_1^{\langle r\rangle})^1_{Q_1} ,\frac12[(W_1^{\langle r\rangle})^0_{Q_1},Q_1Q_2]
+(W_{12}^{\langle r\rangle})^1_{Q_1Q_2}
-\frac12[(W_1^{\langle r\rangle})^0_{Q_1},Q_1Q_2],W_2^{\langle r\rangle})^1_{Q_2}
 \right).
\end{align*}
(Remark: This involves a symmetric ``connection choice'' for $\varkappa^{\langle r\rangle},\covarkappa^{\langle r\rangle}, \coeta^{\langle r\rangle}$.)

Then, we claim, there are equalities
\begin{align*}
\PreAmb_Q^{\langle r\rangle}&:=\im \varkappa^{\langle r\rangle}=\ker\eta^{\langle r\rangle},\\
\Amb_Q^{\langle r\rangle}&:=\im  \lambda^{\langle r\rangle}=\ker \colambda^{\langle r\rangle},\\
\CoAmb_Q^{\langle r\rangle}&:=\im  \covarkappa^{\langle r\rangle}=\ker \varkappa^{\langle r\rangle},\\
\PreCoAmb_Q^{\langle r\rangle}&:=\im \colambda^{\langle r\rangle}=\ker \coeta^{\langle r\rangle}.
\end{align*}

Furthermore, the maps $\varkappa^{\langle r\rangle}$ and $\lambda^{\langle r\rangle}$ induce (inverse) bijections between
$\PreAmb_Q^{\langle r\rangle}$ and $\Amb_Q^{\langle r\rangle}$;
the maps $\colambda^{\langle r\rangle}$ and $\covarkappa^{\langle r\rangle}$ induce (inverse) bijections between
$\CoAmb_Q^{\langle r\rangle}$ and $\PreCoAmb_Q^{\langle r\rangle}$;
leading to a splitting
\[\vFQ^{\langle r\rangle}=\Amb_Q^{\langle r\rangle}\oplus\CoAmb_Q^{\langle r\rangle}\]
with factors
\[\pi_Q^{\langle r\rangle}=  \lambda^{\langle r\rangle}\circ \varkappa^{\langle r\rangle}\qquad
\text{and}\qquad\copi_Q^{\langle r\rangle}=\covarkappa^{\langle r\rangle}\circ \colambda^{\langle r\rangle}.\]

\begin{proof}
The point is that elements  $V^{\langle r\rangle}\in \Amb_Q^{\langle r\rangle}$ are of shape
\[(V^{(r)}/Q)_{\mathrm{split}}
=\begin{bmatrix} 0&0&x_3^{\langle r\rangle}&x_4^{\langle r\rangle}&0&x_6^{\langle r\rangle}&0&x_4^{\langle r\rangle} \end{bmatrix},\]
in bijection to elements $U^{\langle r\rangle}\in \PreAmb_Q^{\langle r\rangle}$
\[U^{\langle r\rangle}=\frac12(x_3^{\langle r\rangle}+x_4^{\langle r\rangle}+x_6^{\langle r\rangle});\]
and elements  $V^{\langle r\rangle}\in \CoAmb_Q^{\langle r\rangle}$ are of shape
\[(V^{(r)}/Q)_{\mathrm{split}}=\begin{bmatrix} x_1^{\langle r\rangle}&x_2^{\langle r\rangle}&0&x_4^{\langle r\rangle}&
x_5^{\langle r\rangle}
&0&x_7^{\langle r\rangle}&-x_4^{\langle r\rangle} \end{bmatrix},\]
in bijection to elements $W^{\langle r\rangle}\in \PreCoAmb_Q^{\langle r\rangle}$
\[W^{\langle r\rangle}=(-2(x_1^{\langle r\rangle}+x_2^{\langle r\rangle}) ,
 (2x_2^{\langle r\rangle}+4x_4^{\langle r\rangle}-2x_7^{\langle r\rangle})Q_1Q_2
 , -2(x_5^{\langle r\rangle}+x_7^{\langle r\rangle})).\]
(But formally, we can just check some relations between
$\eta^{\langle r\rangle},\lambda^{\langle r\rangle},\varkappa^{\langle r\rangle},
\colambda^{\langle r\rangle},\covarkappa^{\langle r\rangle},\coeta^{\langle r\rangle}$.)
\end{proof}
\end{lemma}

\begin{theorem}\label{thm:cpgen}
Suppose that $\Psi$ is an Clifford productive, Clifford conservative, vectorial FQ operation. Then
\begin{equation}
\begin{split}
\Psi(A_1,A_2)=
\ldots(1+U^{\langle3\rangle})(1+U^{\langle2\rangle})(1+U^{\langle1\rangle})\cdot
(Q_1,Q_2)\cdot\\\cdot(1+U^{\langle1\rangle})^{-1}(1+U^{\langle2\rangle})^{-1}(1+U^{\langle3\rangle})^{-1}\ldots,
\end{split}
\label{eq:cpgen}
\end{equation}
where $U^{\langle i\rangle}\in \sFQ^{\langle i\rangle}$.
Furthermore, we can obtain a most economical choice for the $U^{\langle r\rangle}$ (yielding bijective correspondence to $\Psi$)
by choosing $U^{\langle r\rangle}$ from \[\PreAmb_{Q}^{\langle r\rangle}=\{ U^{\langle r\rangle}\in\sFQ^{\langle r\rangle}\,:\,
(U^{\langle r\rangle}){}^{0}_{Q_1}{}^{0}_{Q_2}=0\}.\]
With this choice, $U^{\langle r\rangle}$ can be recovered from $\Psi/\vFQ^{\geq r+1}$ by simple arithmetics.
\begin{proof}
First, on the constructive side:
One can see that with arbitrary choice, \eqref{eq:cpgen} yields a Clifford productive operation $\Psi$.
(Because it makes sense and it behaves so  in every order.)
Let $\Psi^{(r)}$ be the version of the FQ operation when we use only
$U^{\langle2\rangle},\ldots,U^{\langle r\rangle}$, but not further.
Then
\[\Psi^{(r)}(A_1,A_2)-\Psi^{(r-1)}(A_1,A_2)=\Delta\Psi^{(r)} +O((R/Q)^{\geq r+1}),\]
where $\Delta\Psi^{(r)}$ is $r$-homogeneous in $(R/Q)$.
In fact,
\[\Delta\Psi^{(r)}=U^{\langle r\rangle}\cdot(Q_1,Q_2)-(Q_1,Q_2)\cdot U^{\langle r\rangle}.\]
This term describes exactly how much ambiguity arises we step up form order $r-1$ to $r$.
Such ambiguities form $\Amb_Q^{\langle r\rangle}$, and these ambiguities can be achieved
using $U^{\langle r\rangle}\in \PreAmb_Q^{\langle r\rangle} $ (and there is a bijective correspondence).

However, in terms of ambiguities, one cannot do better even in the general case.
Suppose that we an Clifford conservative operation FQ operation from $\vFQ^{\leq r-1}$,
such that it is a Clifford system modulo $\vFQ^{\geq r}$.
Suppose that we manage to extend it to an  FQ operation $\Psi\in \vFQ^{\leq r}$, which is a Clifford system
modulo $\vFQ^{\geq r+1}$.
But there is the possibility of getting another version $\Psi'$ so that
\[\Psi'=\Psi+ V^{\langle r\rangle}+ O((R/Q)^{\geq r+1}),\]
where $V^{\langle r\rangle}\in \vFQ^{\langle r\rangle}$.
Then the Clifford system property implies
\[Q_1V^{\langle r\rangle}_1+V^{\langle r\rangle}_1 Q_1=0,\]
\[Q_1V^{\langle r\rangle}_2+V^{\langle r\rangle}_1 Q_2+Q_2V^{\langle r\rangle}_1+V^{\langle r\rangle}_2 Q_1=0,\]
\[Q_2V^{\langle r\rangle}_2+V^{\langle r\rangle}_2 Q_2=0,\]
which is equivalent to $V^{\langle r\rangle}\in \Amb_{(Q_1,Q_2)}^{\langle r\rangle}$.
Proceeding inductively, we see that every Clifford productive,
Clifford conservative, vectorial FQ operation
occurs in form \eqref{eq:cpgen}, and in fact, using the economical choice for the $U^{\langle r\rangle}$.
Then, in each step, we have a full involutive operation $\Psi^{( r)}$,
and the modifier terms arithmetically determined as follows: If $U^{\langle2\rangle},\ldots,U^{\langle r-1\rangle}$,
are already recovered, then so is $\Psi^{(r-1)}$, and  $U^{\langle r\rangle}=\varkappa_Q^{\langle r\rangle}((\Psi- \Psi^{(r-1)})^{\langle r\rangle})$.
\end{proof}
\end{theorem}
\begin{remark}We could have used the form
\[\begin{split}
\Psi(A_1,A_2)=
(1+U^{\langle1\rangle})(1+U^{\langle2\rangle})(1+U^{\langle3\rangle})\ldots\cdot
(Q_1,Q_2)\cdot\ldots\\\ldots(1+U^{\langle3\rangle})^{-1}(1+U^{\langle2\rangle})^{-1}(1+U^{\langle1\rangle})^{-1};
\end{split}\]
it leads to the same ambiguities at each level.
In this latter form, unicity is even more transparent.
Another possibility is to take the exponential version where  $1+U^{\langle r\rangle}$ is replaced by $\exp U^{\langle r\rangle}$.

\end{remark}
For the sake of the next theorem, we use the general notation $(\Psi_1,\Psi_2) \odot  (\Psi_1,\Psi_2)=
(\Psi_1\Psi_1,\Psi_1\Psi_2+\Psi_2\Psi_1,\Psi_2\Psi_2 ) $.
\begin{theorem} Consider a Clifford conservative, vectorial FQ operation $\Psi$.
Then, we claim, the Clifford productivity of $\Psi$ can be expressed by the equations
\begin{equation}
\copi_Q^{\langle r\rangle}\Psi^{\langle r\rangle}
+\covarkappa_Q^{\langle r\rangle}((\Psi^{\leq r-1}\odot\Psi^{\leq r-1})^{\langle r\rangle})=0
\label{eq:cplay1}
\end{equation}
$(r\geq2)$.
That is, a collection of $\dim \CoAmb_Q^{\langle r\rangle}$ many equations of shape
\begin{equation}
\text{linear combinations of  $p^{[s]}_{\iota_1,\ldots,\iota_r}$'s}=
 \text{nonlinear expressions of  $p^{[s]}_{\varkappa_1,\ldots,\varkappa_h}$'s ($h<r$)}
\label{eq:invlay2} \end{equation}
($r\geq 2$), such that the linear terms on the left are themselves linearly independent formally,
hence  leading to the eliminability of a set of $p^{[s]}_{\iota_1,\ldots,\iota_r}$'s
of cardinality $\dim \CoAmb_Q^{\langle r\rangle}$.
(If
$(\Psi^{\langle r\rangle}/Q)_{\mathrm{split}}=
\begin{bmatrix} \phi_1^{\langle r \rangle}&\phi_2^{\langle r \rangle}&\phi_3^{\langle r \rangle}
&\phi_4^{\langle r \rangle}&\phi_5^{\langle r \rangle}&\phi_6^{\langle r \rangle}&\phi_7^{\langle r \rangle}&\phi_8^{\langle r \rangle} \end{bmatrix},$
then the equations above  are  restrictive equations
for
$  \phi_1^{\langle r \rangle},  \phi_2^{\langle r \rangle}, \phi_4^{\langle r \rangle}- \phi_8^{\langle r \rangle},
\phi_5^{\langle r \rangle},  \phi_7^{\langle r \rangle}.$ )

More precisely, equations \eqref{eq:cplay1} up to $r$ are equivalent to $ (\Psi\odot\Psi)^{\leq r}=(-1,0,-1)$, i.~e.~Clifford
productivity up to order $r$.
\begin{proof}
Indeed, in the light of the proof of the preceding theorem,
involutivity means exactly a determinacy of the coambiguity  terms in
$\vFQ^{\langle r\rangle}= \Amb_Q^{\langle r\rangle} \oplus \CoAmb_Q^{\langle r\rangle}$.
This can be written down by recovering the $U^{\langle i\rangle} $'s, yielding
$\copi_Q^{\langle r\rangle}\Psi^{\langle r\rangle}=\copi_Q^{\langle r\rangle}(\Psi^{(r-1)})^{\langle r\rangle}$, i.~e.
\begin{equation}
\copi_Q^{\langle r\rangle}\Psi^{\langle r\rangle}=\text{something in terms of }\Psi^{\leq r-1}.\label{eq:cplay3}
\end{equation}
However, the exact shape on the right is not particularly straightforward.
On the other hand, the identity $\Psi\odot\Psi=(-1,0,-1)$ implies $(\Psi\odot\Psi)^{\langle r\rangle}=0$,
which implies $ \covarkappa_Q^{\langle r\rangle}(\Psi\odot\Psi)^{\langle r\rangle}=0$, which can be written as \eqref{eq:cplay1}.
Equation \eqref{eq:cplay1} is informative to the same degree as \eqref{eq:cplay3}, yet it cannot add more information to it
(i. e. to Clifford productivity, meant relative to the restrictive equations of lower order);
hence they must be equivalent.
\end{proof}
\end{theorem}
\begin{remark}
(a) In  \eqref{eq:cplay1}, we have used only $\covarkappa_Q^{\langle r\rangle}|_{\PreCoAmb_Q^{\langle r\rangle}}$ essentially.
Indeed, if have the equations up to $r-1$, then
$(\Psi^{\leq r-1}\odot\Psi^{\leq r-1})^{\langle r\rangle}\in\PreCoAmb_Q^{\langle r\rangle}$.
(This is transparent from writing $\Psi^{\leq r-1}=(\Psi^{(r)})^{\leq r}+\Delta\Psi_r$.)

(b) The arguments above have a version in the natural / conjugation-invariant case.
Here one can assume that $(Q_1,Q_2)$ is the Gram-Schmidt orthogonalization of $(A_1,A_2)$ and
$\phi_1^{\langle r\rangle},\ldots,\phi_8^{\langle r\rangle}$ are $r$-homogeneous polynomials of  $r_1,r_2,r_4=-r_8,r_5,r_7$.
(That is, the spaces $\vFQ_Q,\ldots$ are replaced by $\vFQ_{\underline{\mathcal O}^{\mathrm{GS}}},\ldots$)
Everything is completely analogous, except the expansions are in a reduced set of variables  $r_1,r_2,r_4=-r_8,r_5,r_7$.
Or, alternatively, we can assume that $(Q_1,Q_2)$ is the symmetric orthogonalization of $(A_1,A_2)$ and
we deal with polynomials of, say,  $\hat r_1,\hat r_2,\hat r_3,\hat r_4,\hat r_5$; etc.

(c) The setup above is also consistent with symmetry / orthogonal invariance, then $U^{\langle r\rangle}$'s will be
symmetric / orthogonal invariant.
Using the exponential form one can make the setting compatible to transposition invariance.
\end{remark}
\begin{point}

Similar conclusions can be drawn regarding other kinds of Clifford productivity:

If $\Psi$ is a Clifford conservative, Clifford productive, pseudoscalar FQ operation, then
\[\begin{split}
\Psi(A_1,A_2)=
\ldots(1+U^{\langle3\rangle})(1+U^{\langle2\rangle})(1+U^{\langle1\rangle})\cdot
Q_1Q_2\cdot\\\cdot(1+U^{\langle1\rangle})^{-1}(1+U^{\langle2\rangle})^{-1}(1+U^{\langle3\rangle})^{-1}\ldots,
\end{split}\]
with $U^{\langle r\rangle}\in\sFQ^{\langle r\rangle}$   and
 $(U^{\langle r\rangle}){}^{0}_{Q_1 Q_2}=0$.

If $\Psi$ is a  Clifford conservative, floating Clifford productive FQ operation, then
\[\begin{split}
\Psi(A_1,A_2)=
\ldots(1+U_1^{\langle3\rangle})(1+U_1^{\langle2\rangle})(1+U_1^{\langle1\rangle})\cdot
(Q_1,Q_2)\cdot\\\cdot(1+U_2^{\langle1\rangle})^{-1}(1+U_2^{\langle2\rangle})^{-1}(1+U_2^{\langle3\rangle})^{-1}\ldots,
\end{split}\]
with $U_1^{\langle r\rangle},U_2^{\langle r\rangle} \in\sFQ^{\langle r\rangle}$   and
 $(U_1^{\langle r\rangle}){}^{0}_{Q_1 Q_2}=Q_1(U_2^{\langle r\rangle}){}^{0}_{Q_1 Q_2} Q_1^{-1}$.

Again, we can derive the existence of the well-layered restrictive equations; the details are left to the reader.
\end{point}

\section{Involutivity and idempotence}
\begin{point} Typically, we want vectorial FQ operations satisfy one of the
following properties:

\texttt{(Inv)} Involutivity: A vectorial FQ operation $\Xi$ is an involution if
$\Xi\circ\Xi=\Id$.

\texttt{(Idm)} Idempotence :A vectorial FQ operation $\Xi$ is idempotent, if
$\Xi\circ\Xi=\Xi$.

\texttt{(I3)} 3-Idempotence: A vectorial FQ operation $\Xi$ is 3-idempotent, if
$\Xi\circ\Xi\circ\Xi=\Xi$. % (As a last resort, mostly.)

Again, it is reasonable restrict to Clifford conservative operations.
\end{point}
\begin{point}
Suppose that $\Psi$ is a vectorial FQ operation. Using its expansion in the split base, we set
\[\mathrm D\Psi:=\left(
\begin{bmatrix}p^{[1]}_1&p^{[1]}_5\\p^{[2]}_1&p^{[2]}_5\end{bmatrix},
\begin{bmatrix}p^{[1]}_2&p^{[1]}_6\\p^{[2]}_2&p^{[2]}_6\end{bmatrix},
\begin{bmatrix}p^{[1]}_3&p^{[1]}_7\\p^{[2]}_3&p^{[2]}_7\end{bmatrix},
\begin{bmatrix}p^{[1]}_4&p^{[1]}_8\\p^{[2]}_4&p^{[2]}_8\end{bmatrix}
\right).\]
We call this the first differential of $\Psi$.
By direct computation, it is easy to see
\end{point}
\begin{lemma}\label{lem:fidiff}
Suppose that $\Psi_1,\Psi_2$  are  Clifford conservative, vectorial FQ operations. Then
\[\mathrm D(\Psi_2\circ\Psi_1)=\mathrm D\Psi_2.\mathrm D\Psi_1,\]
where the point the right side means multiplication componentwise. \qed
\end{lemma}
In particular, in order to have involutive or idempotent FQ operations in this setting,
the components of $\mathrm D\Psi$ must be involutions or idempotents, respectively.
\begin{cor}
Suppose that $\Psi$ is a transposition invariant and orthogonal invariant natural vectorial FQ operation.

(a) If $\Psi$ is involutive, then $\hat p^{[1]}_1=\hat p^{[1]}_2,\hat p^{[1]}_3,\hat p^{[1]}_4=\hat p^{[1]}_5\in\{-1,1\}$.

(b)  If $\Psi$ is idempotent, then $\hat p^{[1]}_1=\hat p^{[1]}_2,\hat p^{[1]}_3,\hat p^{[1]}_4=\hat p^{[1]}_5\in\{0,1\}$.

(c)  If $\Psi$ is 3-idempotent, then $\hat p^{[1]}_1=\hat p^{[1]}_2,\hat p^{[1]}_3,\hat p^{[1]}_4=\hat p^{[1]}_5\in\{-1,0,1\}$.
\begin{proof}
Direct computation in the mixed base.
\end{proof}
\end{cor}
In these cases we have $8$, $8$ and $27$ choices up to order $1$; which we consider  as the principal types
for involutions, idempotents, and 3-idempotents, respectively.
\begin{point}Consider now expansions around a fixed Clifford system $(Q_1,Q_2)$.
Let $\mathcal E$ be expansion giving the Clifford system itself.
For $L\in\oplus^4\mathrm M_2(\mathbb R)$, let $\mathcal T(L)$ be the expansion which is purely of order $1$ and $\mathrm D(\mathcal  T(L) )=L$.
Then, cf. the previous lemma,
\[(\mathcal E+\mathcal T(L_2))\circ(\mathcal E+\mathcal T(L_1))=(\mathcal E+\mathcal T(L_1L_2)).\]

Let $\Lambda\in\vFQ^{\langle r\rangle}$, such that $r\geq 1$.
Then in the compositions
\[(\mathcal E+\mathcal T(L))\circ(\mathcal E+\Lambda)=\mathcal E+\Glob_L(\Lambda)\]
and
\[(\mathcal E+\Lambda)\circ(\mathcal E+\mathcal T(L))=\mathcal E+\Loc_L(\Lambda),\]
the terms $\Glob_L(\Lambda)$ and $\Loc_L(\Lambda)$ are also from  $\vFQ^{\langle r\rangle}$.
\end{point}
\begin{lemma}
The actions $\Glob_{L_1}$ and $\Loc_{L_2}$ commute with each other.
In fact, there is a two-sided associative action  $\oplus^4\mathrm M_2(\mathbb R)$ on $\vFQ$,
such that   $\Glob_{L_1}$ is a linear right action,
$\Loc_{L_2}$ is a linear left action; moreover, $\Glob$ is also linear in $L_1$.
\begin{proof}
This is the consequence of the associativity of the composition.
\end{proof}
\end{lemma}
\begin{theorem}Suppose that $\Psi$  is a  Clifford conservative, vectorial  FQ operation.
Then $\Psi$ is invertible if and only if $\mathrm D\Psi $ is invertible.
\begin{proof}
$\mathcal E+\mathcal T((\mathrm D\Psi)^{-1})$ will be an inverse modulo $\vFQ^{\geq 2}$.
However, one can easily see that if  $\Phi$ is inverse to $\Psi$  modulo $\vFQ^{\geq r}$,
such that $r\geq2$, then
\[\Phi+ \Loc_{(\mathrm D\Psi)^{-1}} (\Id-\Phi\circ\Psi) \]
is inverse to $\Psi$  modulo $\vFQ^{\geq r+1}$. Taking successive iterations, we obtain an inverse.
\end{proof}
\end{theorem}
\begin{lemma} Suppose that $L\in \oplus^4\mathrm M_2(\mathbb R)$ is an involution.
Let
\[\lambda_L^{\langle r\rangle}:=(\Loc_{L}-\Glob_{L})|_{\vFQ^{\langle r\rangle}},\qquad
\colambda_L^{\langle r\rangle}:=(\Loc_{L}+\Glob_{L})|_{\vFQ^{\langle r\rangle}},\]
\[\varkappa_L^{\langle r\rangle}:=\tfrac14(\Loc_{L}-\Glob_{L})|_{\vFQ^{\langle r\rangle}},\qquad
\covarkappa_L^{\langle r\rangle}:=\tfrac14(\Loc_{L}+\Glob_{L})|_{\vFQ^{\langle r\rangle}}.\]
Then, we claim, there are equalities
\[\Amb_L^{\langle r\rangle}:=\im\,\lambda_L^{\langle r\rangle}=\im\,\varkappa_L^{\langle r\rangle}=
\ker\,\colambda_L^{\langle r\rangle}=\ker\,\covarkappa_L^{\langle r\rangle}\]
and
\[\CoAmb_L^{\langle r\rangle}:=\im\,\colambda_L^{\langle r\rangle}=\im\,\covarkappa_L^{\langle r\rangle}=
\ker\,\lambda_L^{\langle r\rangle}=\ker\,\varkappa_L^{\langle r\rangle}.\]

Moreover, $\lambda_L^{\langle r\rangle}$ and $\varkappa_L^{\langle r\rangle}$
induce inverse bijections on $\Amb_L^{\langle r\rangle}$; and
$\colambda_L^{\langle r\rangle}$ and $\covarkappa_L^{\langle r\rangle}$
induce inverse bijections on $\CoAmb_L^{\langle r\rangle}$.
This leads to a direct decomposition
\[\vFQ^{\langle r\rangle}=\Amb_L^{\langle r\rangle}\oplus\CoAmb_L^{\langle r\rangle}\]
with projection factors
\[\pi_L^{\langle r\rangle}=\varkappa_L^{\langle r\rangle}\circ\lambda_L^{\langle r\rangle}=
\lambda_L^{\langle r\rangle}\circ \varkappa_L^{\langle r\rangle}=
\tfrac12(\id-\Loc_{L}\Glob_{L})|_{\vFQ^{\langle r\rangle}}\]
and
\[\copi_L^{\langle r\rangle} =\covarkappa_L^{\langle r\rangle}\circ\colambda_L^{\langle r\rangle}=
\colambda_L^{\langle r\rangle}\circ \covarkappa_L^{\langle r\rangle}=
\tfrac12(\id+\Loc_{L}\Glob_{L})|_{\vFQ^{\langle r\rangle}}.\]
\begin{proof}
This follows from that $\Glob_L$ and $\Loc_L$ are commuting linear actions with spectrum in $\{-1,1\}$.
\end{proof}
\end{lemma}

\begin{theorem}\label{thm:involgen}
Suppose that $\Psi$ is an involutive, Clifford conservative, vectorial FQ operation with $\mathrm D\Psi=L$. Then
\begin{equation}
\begin{split}
\Psi=\ldots\circ(\Id+U^{\langle4\rangle})\circ(\Id+U^{\langle3\rangle})\circ(\Id+U^{\langle2\rangle})\circ
(\mathcal E+\mathcal T(L))\circ\\
\circ (\Id+U^{\langle2\rangle})^{-1} \circ (\Id+U^{\langle3\rangle})^{-1} \circ(\Id+U^{\langle4\rangle})^{-1}\circ\ldots,
\end{split}
\label{eq:involgen}
\end{equation}
where $U^{\langle i\rangle}\in \vFQ^{\langle i\rangle}$.
Furthermore, we can obtain a most economical choice for the $U^{\langle r\rangle}$ (yielding bijective correspondence to $\Psi$)
by choosing $U^{\langle r\rangle}$   from
$\Amb_L^{\langle r\rangle}$.

With this choice, $U^{\langle r\rangle}$ can be recovered from $\Psi/\vFQ^{\geq r+1}$ by simple arithmetics.

\begin{proof}
First, on the constructive side:
One can see that with arbitrary choice, \eqref{eq:involgen} yields an involutive operation $\Psi$.
(Because it makes sense and it behaves so  in every order.)
Let $\Psi^{(r)}$ be the version of the FQ operation when we use only
$U^{\langle2\rangle},\ldots,U^{\langle r\rangle}$, but not further.
Then
\[\Psi^{(r)}-\Psi^{(r-1)}=\Delta\Psi^{(r)} +O((R/Q)^{\geq r+1}),\]
where $\Delta\Psi^{(r)}$ is $r$-homogeneous in $(R/Q)$. In fact,
\[\Delta\Psi^{(r)}=\Loc_{L}(U^{\langle r\rangle})-\Glob_{L}(U^{\langle r\rangle}).\]
This term describes an ambiguity which arises we step up form order $r-1$ to $r$.
These possible terms form  $\Amb_L^{\langle r\rangle}$.
Furthermore, if an ambiguity $V^{\langle r\rangle}=\Delta\Psi^{(r)}\in\Amb_L^{\langle r\rangle}$ is given, then it is induced by
$U^{\langle r\rangle}=\varkappa_L^{\langle r\rangle}( V^{\langle r\rangle})\in\Amb_L^{\langle r\rangle}  $.
We also see that different choices of $U^{\langle r\rangle}\in\Amb_L^{\langle r\rangle}$ lead to different ambiguities $V^{\langle r\rangle}$.

However, one cannot do better even in the general case.
Suppose that we an Clifford conservative operation FQ operation from $\vFQ^{\leq r-1}$,
with first differential $L$, so its square is in $\Id+\vFQ^{\geq r}$, i. e. it is involutive modulo $\vFQ^{\geq r}$.
Suppose that we manage it extend it to an  FQ operation $\Psi\in \vFQ^{\leq r}$, involutive
modulo $\vFQ^{\geq r+1}$, so its square is in $\Id+\vFQ^{\geq r+1}$.
But there is the possibility of getting another version $\Psi'$ so that
\[\Psi'=\Psi+ V^{\langle r\rangle}+ O((R/Q)^{\geq r+1}),\]
where $V^{\langle r\rangle}\in \vFQ^{\langle r\rangle}$. Then $\Psi'\circ \Psi'\in\Id+\vFQ^{\geq r+1}$ means that
\[\Loc_{L}(V^{\langle r\rangle})+\Glob_{L}(V^{\langle r\rangle})=0,\]
i. e., $V^{\langle r\rangle}\in\Amb_L^{\langle r\rangle}$.
Proceeding inductively, we see that every involutive
Clifford conservative, vectorial FQ operation with first differential $L$
occurs in form \eqref{eq:involgen}, and in fact, using the economical choice for the $U^{\langle r\rangle}$.
Then, in each step, we have a full involutive operation $\Psi^{( r)}$,
and the modifier terms arithmetically determined as follows: If $U^{\langle2\rangle},\ldots,U^{\langle r-1\rangle}$
are already recovered, then $U^{\langle r\rangle}=\varkappa_L^{\langle r\rangle}(\Psi- (\Psi^{(r-1)})^{\langle r\rangle})$.
\end{proof}
\end{theorem}
\begin{cor}
Suppose that $\Psi_1,\Psi_2$  are involutive,  Clifford conservative, vectorial FQ operations. Then,
 $\mathrm D\Psi_1$ and  $\mathrm D\Psi_1$ are conjugates of each other if and only
 $\Psi_1$ and $\Psi_2$ are conjugates of each other by Clifford conservative FQ operations.
\begin{proof}
This is a consequence of the previous theorem and Lemma \ref{lem:fidiff}.
\end{proof}
\end{cor}

\begin{theorem} Consider  vectorial FQ operations $\Psi$, with properties such that $\Psi$ is Clifford conservative
and $\mathrm D\Psi=L$ is an involution.

Then, we claim, the involutivity of $\Psi$ can be expressed by
\begin{equation}
\colambda_L^{\langle r\rangle}\,\Psi^{\langle r\rangle}+\copi_L^{\langle r\rangle}(\Psi^{\leq r-1}\circ\Psi^{\leq r-1})^{\langle r\rangle}=0,
\label{eq:invlay0}
\end{equation}
or
\begin{equation}
\copi_L^{\langle r\rangle}\,\Psi^{\langle r\rangle}
+\covarkappa_L^{\langle r\rangle}(\Psi^{\leq r-1}\circ\Psi^{\leq r-1})^{\langle r\rangle}=0
\label{eq:invlay1}
\end{equation}
$(r\geq2)$.
That is, a collection of $\dim \CoAmb_L^{\langle r\rangle}$ many equations of shape
\begin{equation}
\text{linear combinations of  $p^{[s]}_{\iota_1,\ldots,\iota_r}$'s}=
 \text{nonlinear expressions of  $p^{[s]}_{\lambda_1,\ldots,\lambda_h}$'s ($h<r$)}
\label{eq:invlay2} \end{equation}
($r\geq 2$), such that the linear terms on the left are themselves linearly independent formally
hence  leading to the eliminability of a set of $p^{[s]}_{\iota_1,\ldots,\iota_r}$'s of cardinality $\dim \CoAmb_L^{\langle r\rangle}$.

More precisely, equations \eqref{eq:invlay1}/\eqref{eq:invlay2} up to $r$ are equivalent to $ (\Psi\circ\Psi-\Id)^{\leq r}=0$, i.~e.~involutivity up to order $r$.
\begin{proof}
Indeed, in the light of the proof of the preceding theorem,
involutivity means exactly a determinacy of the coambiguity  terms in
$\vFQ^{\langle r\rangle}= \Amb_L^{\langle r\rangle} \oplus \CoAmb_L^{\langle r\rangle}$.
This can be written down by recovering the $U^{\langle i\rangle} $'s, yielding
$\copi_L^{\langle r\rangle}\Psi^{\langle r\rangle}=\copi_L^{\langle r\rangle}(\Psi^{(r-1)})^{\langle r\rangle}$, i.~e.
\begin{equation}
\copi_L^{\langle r\rangle}\Psi^{\langle r\rangle}=\text{something in terms of }\Psi^{\leq r-1}.\label{eq:invlay3}
\end{equation}
However, the exact shape on the right is not particularly straightforward.
On the other hand, form the identity $\Psi\circ\Psi=\Id$ shows
\[(\Glob_L+\Loc_L)(\Psi^{\langle r\rangle})+(\Psi^{\leq r-1}\circ\Psi^{\leq r-1})^{\langle r\rangle}=0,\]
which implies \eqref{eq:invlay0}, which is equivalent to \eqref{eq:invlay1}.
Equation \eqref{eq:invlay1} is informative to the same degree as \eqref{eq:invlay3}, yet it cannot add more information to it
(i. e. to involutivity);
hence they must be equivalent  (that is relative to the restrictive equations of of lower order).
\end{proof}
\end{theorem}
\begin{remark}
(a) In  \eqref{eq:invlay2}, we have used only $\covarkappa_L^{\langle r\rangle}|_{\CoAmb_L^{\langle r\rangle}}$ essentially.
Indeed, if have the equations up to $r-1$, then
$(\Psi^{\leq r-1}\circ\Psi^{\leq r-1})^{\langle r\rangle}\in\CoAmb_L^{\langle r\rangle}$.
(This is transparent from writing $\Psi^{\leq r-1}=(\Psi^{(r)})^{\leq r}+\Delta\Psi_r$.)

(b) Furthermore, the line of arguments in this section have a version
for conjugation-invariant FQ operations, where the expansions are relative, say, to the
symmetric orthogonalizations, hence  we can assume $\Psi$ has an expansion
in the variables $\hat r_1,\hat r_2,\hat r_3,\hat r_4,\hat r_5$.
However, for the sake of compositions one should consider the conjugation-invariant extensions
$\Psi^{\mathrm{ext}}=\Psi_{\underline{\mathcal O}^{\mathrm{Sy}}}$. In this setting
$\Psi_1\circ\Psi_2$ should be replaced by $\Psi_1^{\mathrm{ext}}\circ\Psi_2$.

Naturally, at the end, this leads to the same result as if we apply  the restrictive equations to more restricted
arguments (whose symmetric orthogonalization is $(Q_1,Q_2)$).

(c) Again, these arguments are compatible to symmetry or orthogonal invariance.
\end{remark}
\begin{example}
Consider the case $\Psi$ is symmetric Clifford conservative,  involutive FQ operation with
\[\hat {\mathbf P}^{[1]}_1=\left[ \begin {array}{cccccccc} -1&-1&-1&-1&-1&1&1&1\end {array}\right].\]

In order to describe  this case, let us consider
the parity twist automorphism $\mathrm{ptw}:\mathfrak C\rightarrow \mathfrak C$ given by
\[\begin {array}{c|cccccccc}
x&1&2&3&4&5&6&7&8\\
\hline
\mathrm{ptw}(x)&2&1&3&4&5&6&8&7
\end {array};\]
and the  indicator function $\mathrm{ind}:\mathfrak C\rightarrow \{1,-1\}$
\[\begin {array}{c|cccccccc}
x&1&2&3&4&5&6&7&8\\
\hline
\mathrm{ind}(x)&-1&-1&-1&-1&-1&1&1&1
\end {array}. \]

The expansion orders $0$ and $1$ are determined, and symmetry implies that it is sufficient to
consider $\hat p_{\iota_1,\ldots,\iota_r}^{[1]}$.
Beyond that, in the mixed base, the involutivity can be expressed by the three sets of equations ($r\geq2$)
\[\hat p_{\iota_1,\ldots,\iota_r}^{[1]}=\text{expression of $\hat p_{\text{lesser that $r$ many indices}}^{[1]}$'s}\tag{I}\]
where $\iota_1*\ldots*\iota_r\in\{4,5\}$, and $\prod_{j=1}^r  \mathrm{ind}(\iota_j)=-1$;
\[\hat p_{\iota_1,\ldots,\iota_r}^{[1]}=\text{expression of $\hat p_{\text{lesser that $r$ many indices}}^{[1]}$'s}\tag{II}\]
where $\{\iota_1,\ldots,\iota_r\}\subset\{3,4,5,6\}$,   $\iota_1*\ldots*\iota_r\in\{3,6\}$, and
$\mathrm{ind}(\iota_1*\ldots*\iota_r)\prod_{j=1}^r  \mathrm{ind}(\iota_j)=1$;
\begin{multline}\hat p_{\iota_1,\ldots,\iota_r}^{[1]}+
\left(\mathrm{ind}(\iota_1*\ldots*\iota_r)
\prod_{j=1}^r  \mathrm{ind}(\iota_j)\right)\hat p_{\mathrm{ptw}(\iota_1),\ldots,\mathrm{ptw}(\iota_r)}^{[1]}=\\=
\text{expression of $\hat p_{\text{lesser that $r$ many indices}}^{[1]}$'s}\tag{III}\end{multline}
where $\{\iota_1,\ldots,\iota_r\}\not\subset\{3,4,5,6\}$,   $\iota_1*\ldots*\iota_r\in\{1,2,3,6,7,8\}$.
In case (III) it can be assumed that $(\iota_1,\ldots,\iota_r)$ is lexicographically precedes
$(\mathrm{ptw}(\iota_1),\ldots,\mathrm{ptw}(\iota_r))$. (It would have also been a good strategy to introduce another base
by ``demixing'' $\{1,2\}$ and $\{7,8\}$.)

In the conjugation-invariant case, we are restricted to $\{\iota_1,\ldots,\iota_r\}\subset\{1,2,3,4,5\}$;
and $\prod_{j=1}^r  \mathrm{ind}(\iota_j)$ simplifies to $(-1)^r$.
\end{example}
\begin{point}
When it comes to idempotent operations, the element $L$ should be idempotent.
In this setting, one should use
\[\lambda^{\langle r\rangle}_L=(\Loc_{L}-\Glob_{L})|_{\vFQ^{\langle r\rangle}},\qquad
\colambda^{\langle r\rangle}_L= (\Loc_{L}+\Glob_{L}-\id)|_{\vFQ^{\langle r\rangle}},\]
\[\varkappa_L^{\langle r\rangle}=(\Loc_{L}-\Glob_{L})|_{\vFQ^{\langle r\rangle}},\qquad
\covarkappa_L^{\langle r\rangle}=(\Loc_{L}+\Glob_{L}-\id)|_{\vFQ^{\langle r\rangle}};\]
with
\[\copi_L^{\langle r\rangle}=(\id-\Loc_{L}-\Glob_{L}+2\Loc_L\Glob_{L})|_{\vFQ^{\langle r\rangle}}.\]
Then we have analogous arguments,
leading to restrictive equations of shape as in \eqref{eq:invlay1}/\eqref{eq:invlay2}.

In the case of 3-idempotents, one can proceed with
\[\lambda^{\langle r\rangle}_L=(\Loc_{L}-\Glob_{L})|_{\vFQ^{\langle r\rangle}},\qquad
\colambda^{\langle r\rangle}_L=(\Loc_{L}^2+\Loc_{L}\Glob_{L}+\Glob_{L}^2-\id)|_{\vFQ^{\langle r\rangle}}, \]
\[\varkappa_L^{\langle r\rangle}=(\Loc_L-\Glob_L+\tfrac34\Loc_L^2\Glob_L-\tfrac34\Loc_L\Glob_L^2)
|_{\vFQ^{\langle r\rangle}},\]
\[\covarkappa_L^{\langle r\rangle}=(
-\tfrac34\Loc_{L}^2\Glob_{L}^2+\tfrac14\Loc_{L}\Glob_{L}+\Glob_{L}^2+\Loc_{L}^2-\id
)|_{\vFQ^{\langle r\rangle}};\]
with
\[\copi_L^{\langle r\rangle}=(\id-
\Loc_L^2-\Glob_L^2+\tfrac32\Loc_L^2\Glob_L^2+\tfrac12\Loc_L\Glob_L)|_{\vFQ^{\langle r\rangle}}.\]

The shape of the restrictive equations is
\[\colambda^{\langle r\rangle}_L\Psi^{\langle r\rangle}+
\copi_L^{\langle r\rangle}(\Psi^{\leq r-1}\circ\Psi^{\leq r-1}\circ\Psi^{\leq r-1})^{\langle r\rangle}=0,\]
or
\[\copi_L^{\langle r\rangle}\Psi^{\langle r\rangle}
+\covarkappa_L^{\langle r\rangle}(\Psi^{\leq r-1}\circ\Psi^{\leq r-1}\circ\Psi^{\leq r-1})^{\langle r\rangle}=0;\]
otherwise, the conclusions are similar.
\end{point}
\begin{cor}
Suppose that $\Psi$ is natural, orthogonal invariant, Clifford conservative vectorial FQ operation.
Then $\Psi$ is Clifford conservative if and only if $\Psi^{\leq1}=(\mathcal O^{\mathrm{Sy}})^{\leq1}$ and $\Psi$ is idempotent.
\begin{proof}
In either way, we have a conjugate of $\mathcal O^{\mathrm{Sy}}$, preserving the noted properties.
\end{proof}
\end{cor}
\begin{example}
Consider the case when $\Psi$ is natural, symmetric Clifford conservative,  Clifford productive FQ operation.
Then, in the mixed base,
\[\hat {\mathbf P}^{[1]}_1=\left[ \begin {array}{cccccccc} 0&0&0&0&0&1&1&1\end {array}\right].\]
The expansion orders $0$ and $1$ are determined,  symmetry implies that it is sufficient to
consider $\hat p_{\iota_1,\ldots,\iota_r}^{[1]}$, and naturality / conjugation invariance implies
that it is sufficient to consider $\{\iota_1,\ldots,\iota_r\}\subset\{1,2,3,4,5\} $.

Beyond that, the Clifford productivity / idempotence  can be expressed by three sets of equations ($r\geq2$)
\[\hat p_{\iota_1,\ldots,\iota_r}^{[1]}=\text{expression of $\hat p_{\text{lesser that $r$ many indices}}^{[1]}$'s}\tag{I}\]
where $\iota_1*\ldots*\iota_r\in\{4,5\}$;
\[\hat p_{\iota_1,\ldots,\iota_r}^{[1]}=\text{expression of $\hat p_{\text{lesser that $r$ many indices}}^{[1]}$'s}\tag{II}\]
where $\{\iota_1,\ldots,\iota_r\}\subset\{3,4,5\}$,   $\iota_1*\ldots*\iota_r\in\{3\}$;
\[\hat p_{\iota_1,\ldots,\iota_r}^{[1]}-
\mathrm{ind}(\iota_1*\ldots*\iota_r)\hat p_{\mathrm{ptw}(\iota_1),\ldots,\mathrm{ptw}(\iota_r)}^{[1]}=
\text{expression of $\hat p_{\text{lesser that $r$ many indices}}^{[1]}$'s}\tag{III}\]
where $\{\iota_1,\ldots,\iota_r\}\not\subset\{3,4,5\}$,   $\iota_1*\ldots*\iota_r\in\{1,2,3,6,7,8\}$.
\end{example}
\section{An example of application of restrictive equations}
\begin{point}
Suppose that we have  some of arithmetical conditions  $\mathcal S$;
and we consider those say, vectorial, FQ operations which satisfy them modulo $\vFQ^{\geq r+1}$.
This yields a subvariety $V_r\subset \vFQ^{\leq r}$.
Taking further restrictions in the expansion we have a sequence of subsets
\[\{*\} \xleftarrow{\theta_{0}}V_r^{\leq 0}\xleftarrow{\theta_1}V_r^{\leq 1}\xleftarrow{\theta_2} \ldots
\xleftarrow{\theta_{r}}V_r^{\leq r}=V_r. \]
We will consider only cases when
\[\tag{AFP}\text{the fibers $\theta_{i}$ are affine linear spaces of constant dimension}\]
of $d_r^{\langle i\rangle}=\dim (\theta_{i})^{-1}(x)$, $x\in V_r^{\leq i-1}$ (affine linear fiber property).
We can collect the dimension data into a table
\[\begin{array}{c|ccccc}
&0&1&2&3&\cdots\\
\hline
0&d_0^{\langle0\rangle}&&&&\\
1&d_1^{\langle0\rangle}&d_1^{\langle1\rangle} &&&\\
2&d_2^{\langle0\rangle}&d_2^{\langle1\rangle}&d_2^{\langle2\rangle}&&\\
3&d_3^{\langle0\rangle}&d_3^{\langle1\rangle}&d_3^{\langle2\rangle}&d_3^{\langle3\rangle}&\\
\vdots&\vdots&\vdots&\vdots&\vdots&\ddots
\end{array}\quad.\]
If $j$ is fixed, then $d_r^{\langle j\rangle}$ is monotone decreasing as $r\rightarrow\infty$.
One can easily see that there are two cases:
The process either runs into inconsistency at some level, or, it converges
\[\lim_{r\rightarrow\infty} d_r^{\langle j\rangle}=d^{\langle j\rangle}_{\mathcal S},\]
to  a sequence of $d^{\langle j\rangle}_{\mathcal S}$'s.
If  the affine linear fiber property holds out, then
one can show that this leads to a possibly
infinite-dimensional filtered variety $V_{\mathcal S}\subset\vFQ$ with fiber dimensions $d^{\langle j\rangle}_{\mathcal S}$.
These numbers tell us the degree of freedom $\mathcal S$ allows on the $j$th expansion level.
The most advantageous case is when $d^{\langle j\rangle}_{\mathcal S}=0$ for all $j$;
this means that  the vectorial operation is completely characterized by $\mathcal S$.

Carrying out a complete analysis along these lines is, in general, is very difficult
(although occasionally works out, cf.~Theorems \ref{thm:biv1} and \ref{thm:biv2});
but doing it up to a finite degree often serves us with useful lessons and ideas.
The following provides an illustration.
\end{point}
\begin{point}
Consider the following situation:
We know the FQ orthogonalization procedure $\mathcal O^{\mathrm{Sy}}$, but we are dissatisfied by it, because its
``axis'' $\mathcal D\circ\mathcal O^{\mathrm{Sy}}$ (a pseudoscalar FQ operation) is not affine scaling invariant.
So, we are looking for a nice FQ orthogonalization procedure $\Psi$ with better properties.

(i) We require that $\Psi$ should be Clifford conservative, transposition invariant, orthogonal invariant, and Clifford productive.
The collection of these properties implies the fiber dimensions
\[\begin{array}{c|cccccc}
&0&1&2&3&4&5\\
\hline
0&0\\
1&0&0\\
2&0&0&2\\
3&0&0&2&12\\
4&0&0&2&12&56\\
5&0&0&2&12&56&270\\
\end{array}\quad.\]
So, it seems, there are many operations like that.
(We know that imposing metric trace commutativity, or the fiber-star property, we would obtain $\mathcal O^{\mathrm{Sy}}$;
but we have chosen an other way.)

(ii) Next, we impose that the pseudoscalar FQ operation $\mathcal D\circ\Psi$
should satisfy the $\hat r_i$-scaling invariances with $i\in\{1,2,3,4,5\}$, with $\alpha=0$ (among them is affine invariance, $i=3,4$).
We have fiber dimensions
\[\begin{array}{c|cccccc}
&0&1&2&3&4&5\\
\hline
0&0\\
1&0&0\\
2&0&0&1\\
3&0&0&0&5\\
4&0&0&0&5&22\\
5&0&0&0&5&16&109\\
\end{array}\quad.\]
We see that the various conditions, well-layered in themselves,
interact with each other, leading to further constraints in \textit{lower} than top expansion degrees (``undercut'').

(iii) Encouraged by the possibilities, next we simply declare that \[\mathcal D\circ\Psi=\underline{\mathcal A}_C,\]
where the central axis operation $\underline{\mathcal A}_C$ is defined by
\begin{equation}
\underline{\mathcal A}_C(A_1,A_2)=\pol\,\frac12(\underline{\mathcal A}_L(A_1,A_2) +\underline{\mathcal A}_R(A_1,A_2)).
\label{eq:axa}\end{equation}
(Assuming that we know from somewhere that $\underline{\mathcal A}_C$ is sufficiently nice for that purpose.)
This leads to fiber dimensions
\[\begin{array}{c|cccccc}
&0&1&2&3&4&5\\
\hline
0&0\\
1&0&0\\
2&0&0&0\\
3&0&0&0&4\\
4&0&0&0&4&16\\
5&0&0&0&4&16&92\\
\end{array}\quad.\]

(iv) Next, we impose the property
\begin{equation}([A_1,\Psi(A_1,A_2)_1]+ [A_2,\Psi(A_1,A_2)_2])^0_{\mathcal A_C(A_1,A_2)}=0, \label{eq:car}\end{equation}
a weakened version of metric trace commutativity.
This leads to fiber dimensions
\[\begin{array}{c|cccccc}
&0&1&2&3&4&5\\
\hline
0&0\\
1&0&0\\
2&0&0&0\\
3&0&0&0&0\\
4&0&0&0&0&0\\
5&0&0&0&0&0&0\\
\end{array}\quad;\]
which is quite encouraging to think that we have characterized a certain FQ operation $\mathcal O^{\mathrm{Sy}, \mathcal A_C}$.
But this is not a proof yet.
However, using \texttt{(CP)} and \eqref{eq:axa}, we can rearrange \eqref{eq:car}, as
\[[(A_1)^1_{\mathcal A_C(A_1,A_2)},\Psi(A_1,A_2)_1]+ [(A_2)^1_{\mathcal A_C(A_1,A_2)},\Psi(A_1,A_2)_2]=0,\]
from which we can deduce
\[\mathcal O^{\mathrm{Sy}, \mathcal A_C}(A_1,A_2)=
\Psi(A_1,A_2)=\mathcal O^{\mathrm{Sy}}( (A_1)^{1}_{\mathcal A_C(A_1,A_2)},(A_2)^{1}_{\mathcal A_C(A_1,A_2)}   ).\]

Here,  the good properties of  $\mathcal O^{\mathrm{Sy}, \mathcal A_C}=\Psi$ defined above are immediate:
Indeed, $\mathcal D\circ \Psi$ will be a Clifford conservative, Clifford productive pseudoscalar operation
multiplicatively commuting with the similar operation $\mathcal A_{C}$, from which one can derive that they are equal.
So, in retrospect, this is quite simple.

(iv') On the other hand, we may look  for an other kind of operation by requiring compatibility with ${\mathcal O}^{\mathrm{afSy}}$ and $\mathcal O^{\mathrm{fSy}}$;
 so, instead of \eqref{eq:car},
\begin{equation}\Psi\circ\mathcal O^{\mathrm{fSy}}=\Psi\qquad
\text{and}
\qquad\Psi\circ{\mathcal O}^{\mathrm{afSy}}=\Psi\label{eq:biaxa}\end{equation}
should be satisfied.
(Unfortunately,  $\mathcal O^{\mathrm{Sy}, \mathcal A_C}$ fails \eqref{eq:biaxa}; and
 $\mathcal O^{\mathrm{Sy}}\circ\mathcal O^{\mathrm{fSy}}$  and $\mathcal O^{\mathrm{Sy}}\circ{\mathcal O}^{\mathrm{afSy}}$ fail \eqref{eq:axa}.)
In this case, the first few fiber dimensions are
\[\begin{array}{c|cccccc}
&0&1&2&3&4&5\\
\hline
0&0\\
1&0&0\\
2&0&0&0\\
3&0&0&0&0\\
4&0&0&0&0&3\\
5&0&0&0&0&3&0\\
\end{array}\quad,\]
leaving the case quite open.
In fact, one can construct a relatively nice operation (axial orthogonalization) satisfying the requirements,
but which we do not explain here.
\end{point}
\begin{point}
At this point, the reader may wonder on the following:

(i) This ``method'' seems to be much more suitable to prove non-existence than existence.
This is true, indeed.
On the other hand, in practice, most results \textit{are} negative.
It is just very easy to make wrong guesses about FQ operations.
In the author's experience, FQ operations almost always fail to be so nice as one would like them to be,
but do not fail to stay a little bit mysterious.
So, in fact, this method serves well our analytical efforts.

(ii)  It is not clear what is the exact relevance of the previous sections in this method, in general.
Indeed, as computations are tedious, sooner or later one should impose non-basic invariance rules, and
one resorts to using computers anyway;  the exact form of the restrictive equations gets irrelevant.
This is also true. On the other hand, in a situation where data grows exponentially in the expansion order,
\textit{any} edge in the computation is welcome.
One should avoid solving large systems of equations as much as possible.
In the author's experience, regarding FQ operations, one cannot really obtain a correct picture
just from extrapolating from a couple of expansion orders (say $r=1,2,3$).
In some relevant cases the first ambiguities appear in expansion orders $r=6,7$.
Hence using appropriate bases, etc., makes a difference.

\end{point}

\end{document}